\documentclass{amsart} % in was article
\usepackage{amssymb}
\usepackage{amsmath}
\usepackage{amsfonts}

\begin{document}
%%%%%%%%%%%%%%%%%%%%%%%%%%%%%%%%%%%%%%%%%%%%%%%%%%%%%%%%%%%%%%%%%%%%%%
%	spaces for your own definitions follows
%%%%%%%%%%%%%%%%%%%%%%%%%%%%%%%%%%%%%%%%%%%%%%%%%%%%%%%%%%%%%%%%%%%%%%
\newtheorem{theo}{Theorem}[section]
\newtheorem{prop}[theo]{Proposition}
\newtheorem{lemma}[theo]{Lemma}
\newtheorem{exam}[theo]{Example}
\newtheorem{coro}[theo]{Corollary}
\theoremstyle{definition}
\newtheorem{defi}[theo]{Definition}
\newtheorem{rem}[theo]{Remark}

%\renewcommand{\theequation}{\mbox{\arabic{section}.\arabic{equation}}}

%letters - added these
\newcommand{\Bb}{{\bf B}}
\newcommand{\Fb}{{\bf F}}
\newcommand{\Nb}{{\bf N}}
\newcommand{\Qb}{{\bf Q}}
\newcommand{\Rb}{{\bf R}}
\newcommand{\Zb}{{\bf Z}}
\newcommand{\Ac}{{\mathcal A}}
\newcommand{\Bc}{{\mathcal B}}
\newcommand{\Cc}{{\mathcal C}}
\newcommand{\Dc}{{\mathcal D}}
\newcommand{\Fc}{{\mathcal F}}
\newcommand{\Ic}{{\mathcal I}}
\newcommand{\Jc}{{\mathcal J}}
\newcommand{\Lc}{{\mathcal L}}
\newcommand{\MC}{{\mathcal M}}
\newcommand{\Oc}{{\mathcal O}}
\newcommand{\Pc}{{\mathcal P}}
\newcommand{\Sc}{{\mathcal S}}
\newcommand{\Tc}{{\mathcal T}}
\newcommand{\Uc}{{\mathcal U}}
\newcommand{\Vc}{{\mathcal V}}
\newcommand{\Xc}{{\mathcal X}}

\newcommand{\ax}{{\rm ax}}
\newcommand{\Acc}{{\rm Acc}}
\newcommand{\Act}{{\rm Act}}
\newcommand{\ded}{{\rm ded}}
\newcommand{\Gm}{{$\Gamma_0$}}
\newcommand{\ID}{{${\rm ID}_1^i(\Oc)$}}
\newcommand{\PAP}{{${\rm PA}(P)$}}
\newcommand{\ACA}{{${\rm ACA}^i$}}
\newcommand{\RefP}{{${\rm Ref}^*({\rm PA}(P))$}}
\newcommand{\RefS}{{${\rm Ref}^*({\rm S}(P))$}}
\newcommand{\Rfn}{{\rm Rfn}}
\newcommand{\tar}{{\rm Tarski}}
\newcommand{\UNFA}{{${\mathcal U}({\rm NFA})$}}

\author{Nik Weaver}

\title [Analysis in $J_2$]
       {Analysis in $J_2$}

\address {Department of Mathematics\\
          Washington University in Saint Louis\\
          Saint Louis, MO 63130}

\email {nweaver@math.wustl.edu}

\date{\em September 12, 2005}

\begin{abstract}
This is an expository paper in which I explain how core mathematics,
particularly abstract analysis, can be developed within a concrete
countable set $J_2$ (the second set in Jensen's constructible
hierarchy). The implication, well-known to proof theorists but probably
not to most mainstream mathematicians, is that ordinary mathematical
practice does not require an enigmatic metaphysical universe of sets.
I go further and argue that $J_2$ is a {\it superior} setting
for normal mathematics because it is free of irrelevant set-theoretic
pathologies and permits stronger formulations of existence results.
\end{abstract}

\maketitle

%%%%%%%%%%%%%%%%%%%%%%%%%%%%%%%%%%%%%%%%%%%%%%%%%%%%%%%%%%%%%%%%%%%%%%%
%	Please insert the article body now
%%%%%%%%%%%%%%%%%%%%%%%%%%%%%%%%%%%%%%%%%%%%%%%%%%%%%%%%%%%%%%%%%%%%%%%

Perhaps many mathematicians would admit to harboring some feelings
of discomfort about the ethereal quality of Cantorian set theory. Yet
draconian alternatives such as intuitionism, which holds that simple
number-theoretic statements like the twin primes conjecture may have
no definite truth value, probably violate the typical working
mathematician's intuition far more severely than any vague unease he
may feel about remote cardinals such as, say, $\aleph_{\aleph_\omega}$.

I believe that ordinary mathematical practice is actually most compatible
with an intermediate foundational stance I call {\it conceptualism}. This
is a modernized version of the {\it predicativist} philosophy originally
formulated by Henri Poincar\'e and Bertrand Russell, which treats
elementary number theory from an essentially platonic point of view but
in the realm of set theory admits only those sets that can be constructed
in a very concrete way. One consequence of this restriction is that on
the conceptualist account all sets are countable. While this would seem
to render the bulk of ordinary mathematics illegitimate, work begun by
Hermann Weyl and continued by many others has gradually shown that in
fact most if not all core mathematics can be developed within surprisingly
minimalist systems. This point has been emphasized in \cite{Fef}.
However, early work of this type tended to be presented
in a somewhat abstruse or idiosyncratic way, while more recent work is
generally couched in axiomatic frameworks that mainstream mathematicians
might find unintuitive. The goal of this paper is to show how to develop
core mathematics in a concrete countable domain called $J_2$ which should
be relatively easily appreciated by mainstream mathematicians with no
special training in logic. (The first two sections of the paper are
probably the main hurdle; however, they might not need to be read in detail
before proceeding to later sections.) Our central novelty is not any major
technical advance but rather the choice of an approach that is closer in
spirit to the classical style, and hence, perhaps, a more congenial
environment for doing normal mathematics than one finds in previous work
of this type (e.g., \cite{CLS, Fef, Grz, Kon, Lor, Sim, Tak2, Wan, Wey, Zah}).
Since a really thorough exposition could fill several volumes I have chosen to
present the foundational material at a fairly high level of detail, followed
by a moderately detailed treatement of basic real analysis and successively
more synoptic treatments of more advanced topics. I focus on abstract analysis
because it tends to have a more set-theoretic flavor than other general areas
of mainstream mathematics, and hence it presents a greater challenge to being
formalized in $J_2$.

Since conceptualism admits very little, if any, of the classical set-theoretic
universe beyond what is used in ordinary mathematics, a case can be made that
not only is it more compatible with the working mathematician's intuition,
it also provides a better fit with actual mathematical practice.
Furthermore, existence
results actually become stronger: in our setting for example, they become
proofs of {\it existence in $J_2$}. An analogy can be drawn with the
transition from the use of naive infinitesimals to the epsilon-delta method,
where the primary motivation for making the switch is to eliminate reliance
on ill-defined metaphysical entities, but a side benefit, apparent only after
one has adopted the new view, is a genuinely deeper and more precise
understanding of the subject matter. For these reasons I believe that
conceptualist systems may come to be seen as a {\it superior} setting for
doing normal mathematics.

I will make some brief comments about the philosophical content of
conceptualism in \S \ref{sect1.4} but otherwise will let the mathematical
development speak for itself. For a more thorough discussion of the ideas
of conceptualism including its philosophical justification, see \cite{Wea1}
(for general readers) or \cite{Wea2} (for readers with some background
in proof theory).

\section{The mini-universe $J_2$}\label{sect1}

\begin{quote}
{\it The main goal of this section is to introduce a countable set $J_2$
which will play the role of a miniature universe in which core mathematics
can be developed. We also formulate a ``first definability principle''
according to which $J_2$ is closed under all normal mathematical
constructions. We briefly discuss some philosophical motivation in
\S \ref{sect1.4}.}
\end{quote}

\subsection{Rudimentary functions}\label{sect1.1}

In order to define $J_2$ we need a way to construct new sets from an
existing repertoire. This is done by means of a class of ``rudimentary''
functions. For more on this material see Chapter VI of \cite{Dev}.

\begin{defi}\label{def1.1} {\it Rudimentary functions} from (tuples of)
sets to sets are constructed according to the following conditions:

(1) The functions
\begin{eqnarray*}
F(x_1, \ldots, x_n) &=& x_i\\
F(x_1, \ldots, x_n) &=& \{x_i, x_j\}\\
F(x_1, \ldots, x_n) &=& x_i - x_j\qquad\hbox{(set-theoretic difference)},
\end{eqnarray*}
where $1 \leq i, j \leq n$, are rudimentary.

(2) If $G$, $H$, and $H_1, \ldots, H_k$ are rudimentary then so are the
functions
\begin{eqnarray*}
F(x_1, \ldots, x_n) &=&
G(H_1(x_1, \ldots, x_n), \ldots, H_k(x_1, \ldots, x_n))\\
F(x_1, \ldots, x_n) &=& \bigcup_{y \in x_1} H(y, x_2, \ldots, x_n).
\end{eqnarray*}

(3) All rudimentary functions are generated by (1) and (2).
\end{defi}

Definitions such as the preceding probably have little meaning the
first time one sees them. Therefore, the reader is invited to build
up his intuition by proving a few cases of the following proposition.

\begin{prop}\label{prop1.2} (\cite{Dev}, Lemma VI.1.1)
The following functions are rudimentary:
\begin{eqnarray*}
F(x) &=& x\\
F(x) &=& \bigcup_{y \in x} y\\
F(x,y) &=& x \cup y\\
F(x_1, \ldots, x_n) &=& \{x_1, \ldots, x_n\}\\
F(x,y) &=& \langle x,y\rangle\quad = \quad\{\{x\},\{x,y\}\}\\
F(x_1, \ldots, x_n) &=& \langle x_1, \ldots, x_n\rangle
\quad = \quad \langle \langle x_1, \ldots, x_{n-1}\rangle, x_n\rangle\\
F(x_1, \ldots, x_n) &=& \{G(y, x_2, \ldots, x_n): y \in x_1\}
\quad\hbox{if $G$ is rudimentary}\\
F(x) &=& \{z \in \bigcup_{y \in x} y: z \in y\hbox{ for all }y \in x\}\\
F(x,y) &=& x \cap y\\
F(x) &=& {\rm 1st}(x) \quad = \quad
\begin{cases}
y& \hbox{if }x = \langle y,z\rangle\hbox{ for some }z\\
\emptyset&\hbox{otherwise}
\end{cases}\\
F(x) &=& {\rm 2nd}(x) \quad = \quad
\begin{cases}
z& \hbox{if }x = \langle y,z\rangle\hbox{ for some }y\\
\emptyset&\hbox{otherwise}
\end{cases}\\
F(x) &=& {\rm dom}(x)\quad = \quad\{y: \langle y,z\rangle \in x
\hbox{ for some }z\}\\
F(x) &=& {\rm ran}(x)\quad = \quad\{z: \langle y,z\rangle \in x
\hbox{ for some }y\}\\
F(x,y) &=& x\times y\\
F(x) &=& x|_y\quad = \quad x \cap (y\times {\rm ran}(x))\\
F(x) &=& \{\langle z,y\rangle: \langle y,z\rangle \in x\}.
\end{eqnarray*}
\end{prop}

Rudimentary functions have an alternative characterization which is
slightly more concrete, if less elegant (\cite{Dev}, Lemma VI.1.11):
a function is rudimentary if and only if it can be built up from the
nine functions
\begin{eqnarray*}
F_0(x,y) &=& \{x,y\}\\
F_1(x,y) &=& x - y\\
F_2(x,y) &=& x \times y\\
F_3(x,y) &=& \{\langle u,z,v\rangle: z \in x\hbox{ and }
\langle u,v\rangle \in y\}\\
F_4(x,y) &=& \{\langle z,v,u\rangle: z \in x\hbox{ and }
\langle u,v\rangle \in y\}\\
F_5(x,y) &=& \{{\rm ran}(x|_z): z \in y\}\\
F_6(x) &=& \bigcup_{y \in x} y\\
F_7(x) &=& {\rm dom}(x)\\
F_8(x) &=& \{\langle u,v\rangle: u, v \in x\hbox{ and }u \in v\}
\end{eqnarray*}
by composition. We will not use this fact, however. In general, we can
avoid dealing with the mechanics of rudimentary functions by invoking
either Proposition \ref{prop1.2} or Corollary \ref{cor1.8} below.

\subsection{The set $J_2$}\label{sect1.2}

We now define the central object of this paper.

\begin{defi}\label{def1.3}
(a) The {\it rudimentary closure} of a set $x$ is the smallest set
$y$ such that $x \subseteq y$, $x \in y$, and $y$ is closed under
application of all rudimentary functions (or equivalently, $y$ is
closed under application of the functions $F_0, \ldots, F_8$
listed above).

(b) $J_0 = \emptyset$; $J_1$ is the rudimentary closure of $J_0$;
$J_2$ is the rudimentary closure of $J_1$.
\end{defi}

These are the first few terms of a transfinite sequence known as
{\it Jensen's hierarchy of constructible sets}. Again, for more on
this material the reader is referred to \cite{Dev}.

$J_1$ has a straightforward description: it is just the set of all
``hereditarily finite'' sets, i.e.,
$$J_1 = \emptyset \cup \Pc(\emptyset) \cup \Pc(\Pc(\emptyset)) \cup\cdots$$
where $\Pc$ denotes the power set operation. Every finite tree gives rise
to a hereditarily finite set by labelling terminal nodes with $\emptyset$
and labelling a general node with the set of labels of its immediate
successors; then the label of the root node (indeed, every label) is
hereditarily finite, and conversely every hereditarily finite set can be
obtained in this way from a finite tree.

We cannot expect to have any simple description of $J_2$, but at least we
have the following basic facts. A set $x$ is said to be {\it transitive}
if $y \in x$ and $z \in y$ implies $z \in x$, i.e., every element of $x$
is a subset of $x$.

\begin{prop}\label{prop1.4}
$J_2$ is countable and transitive.
\end{prop}

Countability follows from the fact that $J_1$ is countable since there are
only countably many rudimentary functions. For transitivity, see Lemma VI.1.7
of \cite{Dev}. The basic idea is to use induction on the minimum number of
times clause (2) of Definition \ref{def1.1} is used in the definition of a
rudimentary function. If we call this the ``degree'' of the function then
we can prove that if a transitive set is closed under application of all
rudimentary functions of degree at most $k$ and we add the value of a
rudimentary function of degree $k+1$ applied to some tuple of elements,
we do not lose transitivity. From this it is possible to inductively infer
that for any $k$ the closure of a transitive set under application of all
rudimentary functions of degree at most $k$ is transitive, which we can then
use to infer transitivity of $J_2$ from transitivity of $J_1$.

\subsection{Definable subsets}\label{sect1.3}

A fundamental property of $J_2$ is the fact that it is closed under the
formation of definable subsets of a set. In order to explain this concept
precisely we need to specify a formal language for set theory. The symbols
of the language we will use are the following:

\begin{quote}
Logical symbols: $\wedge$ (and); $\vee$ (or); $\neg$ (not);
$\Rightarrow$ (implies); $\Leftrightarrow$ (if and only if);
$\forall$ (for all); $\exists$ (there exists).

\noindent Set-theoretic symbols: $\in$; $\subseteq$; $=$.

\noindent Additional symbols: parentheses; a countable list of
variables $v, w, \ldots$.
\end{quote}

This list is redundant in that several symbols could be eliminated
in favor of more complicated expressions involving the other symbols.
But for our purposes there is no particular need to do this.

{\it Atomic formulas} are expressions of the form $v \in w$,
$v \subseteq w$, or $v = w$ for any variables $v$ and $w$. All legal
expressions (formulas of the language) can be built up from the
atomic ones by using the logical symbols in a predictable way.

The interpretation of such an expression is straightforward except
for the following point. In order to interpret a quantifier $\forall v$
or $\exists v$ we need to know the range of possible values of $v$.
Every use of a formula $\phi$ therefore requires that we specify the
intended universe of discourse $x$. We will write $\phi^x$ to explicitly
identify that the range of the variables is $x$, i.e., the possible values
of the quantified variables in $\phi^x$ are precisely the elements of
$x$.  We then say that $\phi$ is {\it relativized to $x$}.

\begin{defi}\label{def1.5}
Let $x$ be a set and let $\phi(v, v_1, \ldots, v_k)$ be a formula
of the language described
above with free (i.e., unquantified) variables $v$, $v_1, \ldots, v_k$.
For each $1 \leq i \leq k$ fix a set $x_i \in x$. Then the set
$$\{y \in x: \phi^x(y, x_1, \ldots, x_k)\}$$
of $y \in x$ which make the expression $\phi$ true when $v_i$
takes the value $x_i$ for $1 \leq i \leq k$ is a {\it definable
subset} of $x$.
\end{defi}

I will try to help the reader develop some intuition for the notion of
definability in \S \ref{sect2.3}.

We have the following result:

\begin{theo}\label{theorem1.6}
If $x \in J_2$ and $y$ is a definable subset of $x$ then $y \in J_2$.
\end{theo}

See Lemma VI.1.17 of \cite{Dev}. This can be proven by induction on
the complexity of the formula $\phi$ which is used to define $y$.
Actually, the right way to prove the theorem is to prove the stronger
statement that if $\phi$ has free variables $w_1, \ldots, w_j,
v_1, \ldots, v_k$ and if $x_1, \ldots, x_k \in x$ then the set
$$\{\langle y_1, \ldots, y_j\rangle \in x^j:
\phi^x(y_1, \ldots, y_j, x_1, \ldots, x_k)\}$$
belongs to $J_2$. One can use specific rudimentary functions
to infer this statement given that it holds for formulas of
lower complexity (in particular, the formulas $\psi_1$ and $\psi_2$
if $\phi = \psi_1 \vee \psi_2$, the formula $\psi$ if $\phi =
(\exists v)\psi$, etc.).

Theorem \ref{theorem1.6} is more powerful than it might
appear. For example, suppose we want to define a subset of $x \in J_2$
using parameters $x_1, \ldots, x_k \in J_2$ some of which might not be
elements of $x$. This can be accomplished using Theorem \ref{theorem1.6}
by replacing $x$ with $x \cup \{x_1, \ldots, x_k\}$, which still belongs
to $J_2$ and makes the $x_i$ available for use as parameters. We may also
need to modify $\phi$ so as to ensure that quantified variables continue to
range only over $x$; this can be done using the {\it bounded quantifiers}
$(\forall v \in x)$ and $(\exists v \in x)$. These are abbreviations,
defined by
$$(\forall v \in x)\Ac \quad\equiv\quad (\forall v)(v \in x \Rightarrow \Ac)$$
and
$$(\exists v \in x)\Ac \quad\equiv\quad (\exists v)(v \in x \wedge \Ac).$$

Even more is true. For example, since $J_2$ is closed under cartesian products
(Proposition \ref{prop1.2}) we can use $x \times y$ in place of $x$ to
obtain definable sets of ordered pairs. Proposition \ref{prop1.2} gives
some indication of the variety of other constructions which are available
in $J_2$. This leads to the following principle.

\begin{defi}\label{def1.7}
A {\it $\Delta_0$ formula} is a formula built up from atomic formulas
using the logical connectives $\wedge, \vee, \neg, \Rightarrow,
\Leftrightarrow$ and the bounded quantifiers $(\forall v \in w)$
and $(\exists v \in w)$ defined above (for any variables $v$ and $w$).
\end{defi}

\begin{coro}\label{cor1.8}
First Definability Principle (FDP): If $x_1, \ldots, x_k, y_1, \ldots,
y_n \in J_2$, $F: J_2^n \to J_2$ is rudimentary, and $\phi$
is a $\Delta_0$ formula with free variables $v, v_1, \ldots, v_k$, then
$$\{z \in F(y_1, \ldots, y_n): \phi(z, x_1, \ldots, x_k)\}$$
belongs to $J_2$.
\end{coro}

(This follows from Lemma VI.1.6 (v) of \cite{Dev}. Note that we need
not relativize $\phi$ because all of its quantifiers are bounded.)

Initially one can only use the FDP by laboriously writing out $\phi$ in
order to ensure that the property one has in mind really can be expressed
as a $\Delta_0$ formula, but it should not take long to develop a good
intuition for which properties are of this type. Very roughly speaking,
a property can be expressed by a $\Delta_0$ formula if one can determine
its truth for given values of the free variables using only sets that
appear as elements of those values. Working through the results of Section
\ref{sect2} should give the reader a basic sense of when this is possible.

\subsection{Motivation}\label{sect1.4}

According to the conceptualist philosophy all sets have to be built up
explicitly from other previously constructed sets. Just what should count
as an ``explicit'' construction may be open to debate, but constructions
using the functions $F_i$ of \S \ref{sect1.1} are certainly acceptable.
Consequently $J_2$ is available for conceptualist mathematics. In contrast,
for example, passing from $\Nb$ to its power set would not be considered
explicit because, although one can understand what subsets of $\Nb$ are,
one has no clear idea of how to concretely generate all of them. For more
on these issues see \cite{Wea1}.

Although there is perhaps no way to establish this decisively, it seems
reasonable to suppose that applying the functions $F_i$, or (what is
equivalent in the long run) passing from a set to its definable subsets,
exhausts the conceptualistically acceptable means of set construction provided
one is allowed to do this transfinitely many times. This claim is similar
in spirit to the Church-Turing thesis which formalizes the informal concept
of an algorithm. One argument for it is empirical since no one has yet
found any other reasonably concrete constructions which cannot be reduced
to these.  Another argument can be made on the basis of the naturality of
alternative characterizations of the definability construction in terms of
infinitary recursion theory \cite{Tak, Tug}.

One's first reaction to the definability construction
may be that it is artificial to exclude sets which cannot be explicitly
defined in some special language. But we are not using an arbitrary language,
we are using the standard language of set theory, and the reason it is the
standard language of set theory is because broadly speaking it is able to
express everything that we can imagine concretely doing with sets. Thus,
the possibility that there might exist ``impredicative'' sets that one in
principle cannot imagine actually constructing hinges on
the philosophical question of whether there is some sense in which the
universe of sets can be conceived of as a well-defined, independently
existing, canonical entity to which we have, even in principle, only
limited access.

In any case, a conceptualist should be able to work not only with $J_2$,
but also with $J_\omega$, $J_{\omega^2}$, $J_{\omega^\omega}$, etc. Just
how far into the transfinite he should be willing to go is a rather
subtle question that I discuss in detail in \cite{Wea2}. However, from
the point of view of ordinary mainstream mathematics the question is not
pressing since for these purposes virtually everything one can do in
any $J_\alpha$ can already be done in $J_2$. But it is good to keep in
mind the possibility of going further should the need arise. For instance,
this might be necessary in order to carry out some transfinite construction
which is not covered by the FDP.

Some readers may be more familiar with G\"odel's hierarchy $L_\alpha$
of constructible sets, so I should mention that we could have used these
instead (with $\alpha$ a limit ordinal). The $J_\alpha$ are merely a
little more convenient since they have better closure properties.

\section{Set Theory}\label{sect2}

\begin{quote}
{\it From here on we will be working within $J_2$. We begin in this
section by defining notions of ``sets'' and ``classes'', which we call
``$\iota$-sets'' and ``$\iota$-classes'', that are appropriate to this
context. We formulate a second definability principle relevant to
$\iota$-classes and define and establish basic facts about $\iota$-functions
and $\iota$-relations. In \S \ref{sect2.3} and \S \ref{sect2.4} we observe
that in the world of $J_2$ all sets are countable and there exists
a universal well-ordering.}
\end{quote}

\subsection{Sets and classes}\label{sect2.1}

According to \S \ref{sect1.4}, conceptualistically the set-theoretic universe
may be conceived as being built up in stages $J_\alpha$ indexed by ordinals
$\alpha$. Thus, if we stop the construction at any $\alpha$, our notion
of ``set'' at that point will coincide with the elements of $J_\alpha$.
Just as in classical set theory, it is also convenient to have a notion
of ``classes'' which may or may not be sets. It might seem natural in this
setting to take as classes all definable subsets of the current universe.
However, this is not a good definition because, for example, $\Nb$ can
realized as a set in $J_2$ (see \S \ref{sect2.3}) and the fact that $J_2$
is countable permits a diagonalization argument by which we can construct
a subset of $\Nb$ which is a definable subset of $J_2$ but not an element
of $J_2$. If we do not want proper classes to possibly be contained in sets,
the right definition seems to be the following.

\begin{defi}\label{def2.1}
An {\it $\iota$-set} is an element of $J_2$. An {\it $\iota$-subset} is
a subset that is an $\iota$-set. An {\it $\iota$-class} is a definable
subset of $J_2$ whose intersection with every $\iota$-set is an
$\iota$-set, and a {\it proper $\iota$-class} is an $\iota$-class that
is not an $\iota$-set. An {\it $\iota$-subclass} is a subset that is an
$\iota$-class.
\end{defi}

The requirement that the intersection of an $\iota$-class with any $\iota$-set
must be an $\iota$-set imposes a kind of uniformity condition on
$\iota$-classes. The presence of extra uniformity conditions will be
a recurring motif throughout this paper.

\begin{prop}\label{prop2.2}
The following hold.

(a) Every $\iota$-set is an $\iota$-class, every element of an $\iota$-class
is an $\iota$-set, and every $\iota$-class contained in an $\iota$-set is an
$\iota$-set.

(b) Let $x$ and $y$ be $\iota$-sets. Then $F(x,y)$ is an $\iota$-set for any
rudimentary function $F$; in particular, $x \cup y$, $x \cap y$, $x - y$,
$x \times y$, $\bigcup_{z \in x} z$, and (provided $x \neq \emptyset$)
$\bigcap_{z \in x} z$ are $\iota$-sets.

(c) Let $X$ and $Y$ be $\iota$-classes. Then $X \cup Y$, $X \cap Y$,
$J_2 - X$, and $X \times Y$ are $\iota$-classes.

(d) If $X$ is an $\iota$-class then so is
$$\Pc_\iota(X) = \{x \in J_2: x \subseteq X\}.$$

(e) If $x$ is an $\iota$-set and $Y \subseteq x \times J_2$ is an
$\iota$-class then $\bigcup Y_a$ and $\bigcap Y_a$ are $\iota$-classes,
where $Y_a = \{y:\langle a,y\rangle \in Y\}$ for all $a \in x$.

(f) If $x$ is an $\iota$-set and $Y \subseteq x \times J_2$ is an
$\iota$-class then
\begin{eqnarray*}
&{\prod}_\iota Y_a = \{f \in J_2:
f\hbox{ is a function with domain $x$}&\\
&\qquad\qquad\qquad\qquad\qquad\qquad
\hbox{such that $f(a) \in Y_a$ for all }
a \in x\}&
\end{eqnarray*}
is an $\iota$-class.
\end{prop}

\begin{proof}
(a) The first statement follows from transitivity of $J_2$
(Proposition \ref{prop1.4}) since for any $x \in J_2$ we have
$$x = \{y \in J_2: y \in x\},$$
which shows that $x$ is a definable subset of $J_2$. Since the
intersection of any two $\iota$-sets is an $\iota$-set (see part (b)),
this implies that $x$ is an $\iota$-class. The second and third
assertions are trivial.

(b) By definition, $J_2$ is closed under application of all rudimentary
functions, so this follows immediately from the definition of $\iota$-sets.

(c) Let $X$ and $Y$ be $\iota$-classes. We will first show that $X \times Y
\subseteq J_2$ is definable; the other cases are similar but easier.

Since $X$ and $Y$ are $\iota$-classes, we have
$$X = \{x \in J_2: \phi^{J_2}(x,x_1, \ldots, x_k)\}$$
and
$$Y = \{y \in J_2: \psi^{J_2}(y,y_1, \ldots, y_l)\}$$
for some $x_1, \ldots, x_k, y_1, \ldots, y_l \in J_2$ and some formulas
$\phi$ and $\psi$. Then
$$X \times Y =
\{z \in J_2: (\exists x, y \in J_2)(z = \langle x,y\rangle
\wedge \phi^{J_2}(x, x_1, \ldots, x_k)
\wedge \psi^{J_2}(y, y_1, \ldots, y_l))\}.$$
This works because every ordered pair in $X \times Y$ belongs to
$J_2$ by the fact that the function $F(x,y) = \langle x,y\rangle$
is rudimentary (Proposition \ref{prop1.2}).

In the above formula we used the statement ``$z = \langle x,y\rangle$''.
We must verify that this could be expressed in the language specified
in \S \ref{sect1.3}. It can be written out as follows:
\begin{eqnarray*}
&(\exists v)(\exists w)\big((\forall u)[u \in v \Leftrightarrow u = x]
\wedge (\forall u)[u \in w \Leftrightarrow (u = x \vee u = y)]&\\
&\qquad\qquad\qquad\qquad\qquad\qquad\qquad\qquad\qquad
\wedge (\forall u)[u \in z \Leftrightarrow (u = v \vee u = w)]\big).&
\end{eqnarray*}
This asserts that $z = \{v, w\}$ where $v = \{x\}$ and
$w = \{x, y\}$.

We omit the simple proofs that $X \cup Y$, $X \cap Y$, and $J_2 - X$ are
definable. The fact that the intersection of each of these with any
$\iota$-set is an $\iota$-set follows easily from part (b) together with
the fact that this is true of $X$ and $Y$. So each is an $\iota$-class.

To see that $X \times Y$ is an $\iota$-class, consider its intersection
with an $\iota$-set $x$. By Proposition \ref{prop1.2}, ${\rm dom}(x)$ and
${\rm ran}(x)$ are $\iota$-sets, and we have
$$(X \times Y) \cap x
= ((X \cap {\rm dom}(x))\times (Y \cap {\rm ran}(x))) \cap x.$$
Since $X$ and $Y$ are $\iota$-classes, part (b) now shows that
the final expression is an $\iota$-set, as desired.

(d) To see that $\Pc_\iota(X) \subseteq J_2$ is definable, say
$X = \{x \in J_2: \phi^{J_2}(x,x_1, \ldots, x_k)\}$
for some $x_1, \ldots, x_k \in J_2$ and some formula $\phi$. Then
$$\Pc_\iota(X) = \{y \in J_2: (\forall x \in J_2)(x \in y \Rightarrow
\phi^{J_2}(x,x_1, \ldots, x_k))\}.$$

Now let $x$ be an $\iota$-set. Then $X \cap \bigcup_{y \in x} y$ is an
$\iota$-set by part (b) and hence
$$\Pc_\iota(X) \cap x = \{z \in x: z \subseteq X\}
= \{z \in x: z \subseteq X \cap \bigcup_{y \in x} y\}$$
is an $\iota$-set by the FDP (Corollary \ref{cor1.8}), so $\Pc_\iota(X)$ is an
$\iota$-class.

(e) The proof that $\bigcup Y_a$ and $\bigcap Y_a$ are definable is similar
to the proof of the corresponding fact for $X \times Y$ given in part (c). To
show that $\bigcup Y_a$ is an $\iota$-class, let $z$ be an $\iota$-set; then
$$z \cap \bigcup Y_a = {\rm ran}(Y \cap (x \times z)),$$
which is an $\iota$-set by part (b) and the fact that $Y$ is an $\iota$-class.
Applying this result to the $\iota$-class $Z = (x \times J_2) - Y$ yields
that $\bigcap Y_a$ is an $\iota$-class since $z \cap \bigcap Y_a =
z - (z \cap \bigcup Z_a)$.

(f) The proof of definability of ${\prod}_\iota Y_a$ is again similar to the
corresponding proof in part (c). Now consider its intersection with an
$\iota$-set $z$. Let
$$y = Y \cap \left(x \times {\rm ran}\left(\bigcup_{a \in z} a\right)\right);$$
this is an $\iota$-set by part (b) of this proposition and the fact that
$Y$ is an $\iota$-class. Then
$$z \cap {\prod}_\iota Y_a = \{f \in z:
f\hbox{ is a function with domain $x$ such that }f \subseteq y\}.$$
Thus, to show that $z \cap {\prod}_\iota Y_a$ is an $\iota$-set it will
suffice to show that $\{f \in z: f$ is a function with domain $x\}$ is an
$\iota$-set, since its intersection with $\Pc_\iota(y)$ will then also be
an $\iota$-set by the fact that $\Pc_\iota(y)$ is an $\iota$-class (part (d)).
We can do this using the FDP because $\{f \in z: f$ is a function with domain
$x\}$ equals
$$\left\{f \in z:
f \subseteq x \times {\rm ran}\left(\bigcup_{y \in z} y\right) \wedge
(\forall v \in x)\left(\exists ! w \in {\rm ran}\left(\bigcup_{y \in z} y\right)\right)(\langle v,w\rangle \in f)\right\}.$$
Here $(\exists! w \in v) \Ac(w)$ abbreviates
$$(\exists w \in v)(\Ac(w) \wedge (\forall a \in v)(\Ac(a) \Rightarrow
a = w))$$
and we can express $\langle v,w\rangle \in f$ in $\Delta_0$ form as
\begin{eqnarray*}
&(\exists b \in f)(\exists c,d \in b)[v \in c \wedge v \in d
\wedge w \in d \wedge (\forall e \in b)(e = c \vee e = d)&\\
&\qquad\qquad\qquad\qquad\qquad\qquad\qquad
\wedge (\forall e \in c)(e = v)
\wedge (\forall e \in d)(e = v \vee e = w)].&
\end{eqnarray*}
This completes the proof.
\end{proof}

In general, verifying that a putative $\iota$-class is a definable subset of
$J_2$ is usually simply a matter of writing out its description in the formal
language specified in \S \ref{sect1.3}. The only common complication is the
possible need to use ordered pairs, which can be handled as in the proof of
Proposition \ref{prop2.2} (c) above. So from here on I will generally
simply say ``verifying definability is straightforward'' when this
conclusion is needed.

Checking that the intersection of a putative $\iota$-class with any
$\iota$-set is an $\iota$-set is not always quite so simple. Usually
the quickest way to do this is to intersect with an arbitrary $\iota$-set
and apply the FDP. We also have the following general tool which is
sometimes useful.

\begin{prop}\label{prop2.3}
Second Definability Principle (SDP): Let $X$ be an $\iota$-class and let $Y$
be a subset of $X$ which is definable by means of a $\Delta_0$ formula with
parameters from $J_2$. Then $Y$ is an $\iota$-class.
\end{prop}

If $X = J_2$ then this is a straightforward consequence of the FDP which
applies that principle to $Y \cap x$ for arbitrary $x \in J_2$; see
(\cite{Dev}, Lemma VI.1.6 (v)). In the general case we have $Y = X \cap Y'$
where $Y'$ is the set which results from replacing $X$ with $J_2$ in the
definition of $Y$. So we reduce to the previous case using Proposition
\ref{prop2.2} (c).

\subsection{Functions and relations}\label{sect2.2}

Along with $\iota$-sets and $\iota$-classes we have analogous notions of
$\iota$-relations and $\iota$-functions. I use the notations
$F[X] = \{F(a): a \in X\}$ and $F^{-1}[Y] = \{a \in X: F(a) \in Y\}$.
The {\it graph} of a function $F: X \to Y$ is the set
$$\Gamma(F) = \{\langle a,b\rangle \in X \times Y: b = F(a)\}$$
(literally the same thing as $F$, but I will follow ordinary
mathematical usage which distinguishes them).

\begin{defi}\label{def2.4}
An {\it $\iota$-relation} on an $\iota$-class $X$ is an $\iota$-class
contained in $X \times X$. An {\it $\iota$-function} from an $\iota$-class
$X$ to an $\iota$-class $Y$ is a function $F: X \to Y$ such that
(1) $\Gamma(F)$ is an $\iota$-class and (2) for each $\iota$-subset
$x \subseteq X$, $F[x]$ is contained in an $\iota$-set. An
{\it $\iota$-injection (surjection)} is an $\iota$-function which is
injective (surjective). An {\it $\iota$-bijection} is a bijection which
is an $\iota$-function in both directions.

A {\it small $\iota$-relation} is an $\iota$-relation which is an $\iota$-set
and a {\it small $\iota$-function} is an $\iota$-function whose graph is an
$\iota$-set.
\end{defi}

We can usually verify condition (2) on $\iota$-functions using the FDP,
since $F[x]$ will in fact be an $\iota$-set (Proposition \ref{prop2.6} (b)).

We collect basic properties of $\iota$-relations and $\iota$-functions in the
following two propositions.

\begin{prop}\label{prop2.5}
Let $X$ be an $\iota$-class and let $R$ be an $\iota$-relation on $X$.

(a) $R$ is small if and only if ${\rm dom}(R)$ and ${\rm ran}(R)$ are
$\iota$-sets.

(b) $R^{-1} = \{\langle b,a\rangle: \langle a,b\rangle \in R\}$ is also
an $\iota$-relation on $X$.

(c) If $Y \subseteq X$ is an $\iota$-class then $R \cap (Y \times Y)$ is an
$\iota$-relation on $Y$.
\end{prop}

\begin{proof}
(a) If $R$ is small then its domain and range are $\iota$-sets by Proposition
\ref{prop2.2} (b) since the domain and range functions are rudimentary
(Proposition \ref{prop1.2}). Conversely, if its domain and range are
$\iota$-sets then their product is an $\iota$-set by Proposition
\ref{prop2.2} (b) and $R$ is then an $\iota$-set by Proposition
\ref{prop2.2} (a).

(b) Definability of $R^{-1}$ is straightforward. Now let
$F(x) = \{\langle b,a\rangle: \langle a,b\rangle \in x\}$; by Proposition
\ref{prop1.2} this is a rudimentary function. Fix an $\iota$-set $x$. Then
$R^{-1} \cap x = F(R \cap F(x))$. Now $F$ takes $\iota$-sets to $\iota$-sets
because it is rudimentary, and $R$ is an $\iota$-class, so this shows that
$R^{-1} \cap x$ is an $\iota$-set. We conclude that $R^{-1}$ is an
$\iota$-class.

(c) This follows from the intersection and product clauses of Proposition
\ref{prop2.2} (c).
\end{proof}

\begin{prop}\label{prop2.6}
Let $X$ and $Y$ be $\iota$-classes and let $F: X \to Y$ be an $\iota$-function.

(a) $F$ is a small $\iota$-function if and only if $X$ is an $\iota$-set.

(b) If $X_0 \subseteq X$ is an $\iota$-class then $F|_{X_0}$ is an
$\iota$-function. If $x \subseteq X$ is an $\iota$-set then $F[x]$
is an $\iota$-set.

(c) If $Y_0 \subseteq Y$ is an $\iota$-class then $F^{-1}[Y_0]$ is
an $\iota$-class. If $F$ is a small $\iota$-function and $y \subseteq Y$
is an $\iota$-set then $F^{-1}[y]$ is an $\iota$-set.

(d) If $Z$ is an $\iota$-class and $G: Y \to Z$ is an $\iota$-function then
$G\circ F: X \to Z$ is an $\iota$-function.

(e) If $Z$ is an $\iota$-class then the identity map from $Z$ to itself is a
$\iota$-function, as is any constant function on $Z$.
\end{prop}

\begin{proof}
(a) If $F$ is a small $\iota$-function then $X = {\rm dom}(\Gamma(F))$ is
an $\iota$-set since the domain function is rudimentary (Proposition
\ref{prop1.2}). Conversely, if $X$ is an $\iota$-set then we can find an
$\iota$-set $y$ such that $F[X] \subseteq y$; then $\Gamma(F) \subseteq
X \times y$ is an $\iota$-set by Proposition \ref{prop2.2} (a) and (b).

(b) $F|_{X_0}$ trivially satisfies condition (2) for $\iota$-functions given
that $F$ does, and its graph is an $\iota$-class because $\Gamma(F|_{X_0}) =
\Gamma(F) \cap (X_0 \times J_2)$ is an $\iota$-class by Proposition
\ref{prop2.2} (c). It follows from this and part (a) that if $x \subseteq X$
is an $\iota$-set then $F|_x$ is a small $\iota$-function; then since
$F[x] = {\rm ran}(\Gamma(F|_x))$ and the range function is rudimentary
(Proposition \ref{prop1.2}) we infer that $F[x]$ is an $\iota$-set.

(c) We prove the second statement first. If $F$ is a small $\iota$-function
then $X$ is an $\iota$-set by part (a); if $y \subseteq Y$ is also an
$\iota$-set then $\Gamma(F) \cap (X \times y)$ is an $\iota$-set with
domain $F^{-1}[y]$. This is then an $\iota$-set since the domain function
is rudimentary.

Definability of $F^{-1}[Y_0]$ is straightforward. Now consider its
intersection with an $\iota$-set $x$. Without loss of generality assume
$x \subseteq X$ and find an $\iota$-set $y$ such that $F[x] \subseteq y$.
Then
$$F^{-1}[Y_0] \cap x = F|_x^{-1}[Y_0 \cap y];$$
since $F|_x$ is a small $\iota$-function (by parts (a) and (b)) and
$Y_0 \cap y$ is an $\iota$-set, this is an $\iota$-set by what we just
proved.

(d) First suppose $F$ is a small $\iota$-function. Then $X$, $y = F[X]$, and
$z = G[y]$ are $\iota$-sets and $G|_y$ is a small $\iota$-function
by parts (a) and (b), and this implies that $\Gamma(G \circ F)$ is an
$\iota$-set by the FDP since
\begin{eqnarray*}
\Gamma(G \circ F) &=& \big\{a \in X \times z:\\
&&\!\!\!\!\!\!\!\!\!\!\!\!\!\!\!\!\!
(\exists u \in X)(\exists v \in y)
(\exists w \in z)(\langle u,v\rangle \in \Gamma(F)
\wedge \langle v,w\rangle \in \Gamma(G|_y)
\wedge a = \langle u,w\rangle)\big\}.
\end{eqnarray*}
Here ordered pairs are handled as in the proof of Proposition \ref{prop2.2}
(f).

Now suppose $F$ is an $\iota$-function. Definability of $\Gamma(G \circ F)$
is shown by a straightforward modification of the above expression and
condition (2) on $\iota$-functions is trivial. To verify that
$\Gamma(G \circ F)$ is an $\iota$-class, fix an $\iota$-set $x$; then
$$\Gamma(G \circ F) \cap x =
\Gamma(G \circ F|_{{\rm dom}(x)}) \cap x.$$
But $F|_{{\rm dom}(x)}$ is a small $\iota$-function since
the domain function is rudimentary (Proposition \ref{prop1.2}), so
$\Gamma(G\circ F|_{{\rm dom}(x)})$ is an $\iota$-set by what we proved first.
Hence $\Gamma(G\circ F) \cap x$ is an $\iota$-set.

(e) Let $F$ be the identity map on $Z$. It is straightforward to check that
$\Gamma(F) = \{\langle a,a\rangle: a \in Z\}$ is definable, and $F$ trivially
satisfies condition (2) on $\iota$-functions. To see that $\Gamma(F)$ is an
$\iota$-class, let $x$ be an $\iota$-set. Then
$$\Gamma(F) \cap x = \{\langle a,a\rangle: a \in Z \cap {\rm dom}(x)\}$$
is an $\iota$-set, as desired. Any constant function $F$ whose image is
an $\iota$-set $a$ has graph $Z \times \{a\}$, which is an $\iota$-class by
Proposition 2.2 (c), and it easily follows that $F$ is an $\iota$-function.
\end{proof}

\begin{prop}\label{prop2.7}
Every rudimentary function is an $\iota$-function.
\end{prop}

This follows from (\cite{Dev}, Lemma VI.1.3).

\subsection{Natural numbers and countability}\label{sect2.3}

We encode the natural numbers as sets in the standard way by inductively
setting
$$0 = \emptyset\qquad{\rm and}\qquad n+1 = \{0,\ldots, n\}.$$

\begin{prop}\label{prop2.8}
$\Nb = \{0, 1, 2, \ldots\} = \{\emptyset, \{\emptyset\},
\{\emptyset, \{\emptyset\}\}, \ldots, \}$ is an $\iota$-set. The standard
operations $+$ and $\cdot$ on $\Nb$ are small $\iota$-functions and the
standard ordering $\leq$ on $\Nb$ is a small $\iota$-relation.
\end{prop}

\begin{proof}
The first claim can be proven using Theorem \ref{theorem1.6} by
showing that $\Nb$ is a definable subset of $J_1$. It is clear that
$\Nb \subseteq J_1$ from the characterization of $J_1$ as the set of all
hereditarily finite sets. Moreover, every $n \in \Nb$ is transitive
($x \in n$ and $y \in x$ implies $y \in n$) and linearly ordered by
$\in$ (for any $x,y \in n$, either $x \in y$ or $y \in x$ or $x = y$),
and it is not hard to see that any hereditarily finite set with these
properties belongs to $\Nb$. (Start by observing that the least element
of such a set must be $\emptyset$.) Thus one can write a formula $\phi$
that expresses transitivity and linearity and obtain
$\Nb = \{x \in J_1: \phi^{J_1}(x)\}$, and hence $\Nb$ is an $\iota$-set
by Theorem \ref{theorem1.6}.

The standard ordering $\leq$ on $\Nb$ actually coincides with $\subseteq$,
so one sees that $\leq$ is an $\iota$-set by writing
$$\leq\, = \{x \in \Nb \times \Nb:
(\exists m, n \in \Nb)(x = \langle m,n\rangle \wedge m \subseteq n)\}$$
and using the FDP. (As usual, the ordered pair notation is an abbreviation;
see the proof of Proposition \ref{prop2.2} (f).) Realizing the graphs of
$+$ and $\cdot$ as $\iota$-sets is a little more complicated (but this
immediately implies that they are $\iota$-functions); I will describe a
method for getting the graph of $+$ using the FDP. First, the defining
formula for the graph of $+$ should assert that $x = \langle a,b,c\rangle$
is an ordered triple of elements of $\Nb$. Next, it should assert that if
$b = 0$ then $a = c$. Finally, there should be an inductive clause that
handles the case $b > 0$ which should state ``if $b \neq 0$ then there exists a
function $f \in J_1$ with domain $a + 1$ ($= a \cup \{a\}$) such that
$f(0) = b$,
$f(a) = c$, and for any $n < a$ we have $f(n+1) = f(n) \cup \{f(n)\}$.''
The existence of such a function ensures that $a + b = c$, and the point
of taking $f \in J_1$ is that $J_1$ contains all finite partial functions
from $\Nb$ to $\Nb$, so that if $a + b = c$ then there is such a
function $f$ in $J_1$. This shows that $+$ is a small $\iota$-function,
and the argument for $\cdot$ is analogous; its defining formula can be
built up using the graph of $+$ as a parameter.
\end{proof}

Now that the familiar operations of arithmetic are available in $J_2$, it
is easier to get some idea of the scope of the concept of definability.
For example, the set of even numbers is defined by the expression
$$\{v \in \Nb: (\exists w \in \Nb)(2\cdot w = v)\},$$
the set of prime numbers is defined by the expression
$$\left\{u \in \Nb: \neg(u = 1) \wedge
(\forall v, w \in \Nb)\left[v\cdot w = u \Rightarrow
(v = 1 \vee w = 1)\right]\right\},$$
and so on. Of course any finite set of natural numbers is definable by
the formula
$$\{v \in \Nb: (v = \_\_) \vee \cdots \vee (v = \_\_)\}$$
with arbitrary numbers inserted in the blank spots. In fact, practically
all sets of natural numbers that we can in any sense explicitly describe
are in $J_2$. But not all: by enumerating all $x \in J_2$ contained
in $\Nb$ (recall that $J_2$ is countable) and diagonalizing we can get
$X \subseteq \Nb$ with $X \in J_3$ but $X \not\in J_2$, and we can go on
to get subsets of $\Nb$ in $J_4$ but not $J_3$, etc. Such sets can be
specified precisely and in some sense explicitly but they are rather
unusual and not likely to occur in ordinary mathematical practice.

One might be concerned that even if we rarely if ever need to explicitly
specify such an apparently pathological object as a set of numbers that
is not in $J_2$, these might still arise as solutions to problems of
interest. We must ask whether any standard existence proofs fail in
$J_2$ because the solutions which they nonconstructively identify might
not belong to $J_2$. The answer to this question is that in general
even nonconstructive arguments will not lead one out of $J_2$; this is not
an absolute fact but rather an empirical observation about the kinds of
arguments used in ordinary mathematics. In other words, even nonconstructive
existence proofs generally {\it do} identify definable solutions to the
problems they address. This is the meaning of the comment I made in the
introduction about existence results being strengthened when one does
mathematics in $J_2$.

We now want to observe that within $J_2$ every $\iota$-set is countable,
i.e., in bijection either with a natural number or with $\Nb$. We already know
this is true ``from the outside'' as it were by Proposition \ref{prop1.4}, but
the new claim is that for every $\iota$-set in $J_2$ there is such a
bijection {\it within $J_2$}, i.e., an $\iota$-bijection. A key tool in
the proof is the following result:

\begin{lemma}\label{lemma2.9}
Suppose $x,y \in J_2$, $x$ is an infinite subset of $y$, and there is a
surjective $\iota$-function $f$ from $\Nb$ onto $y$. Then there is an
$\iota$-bijection from $\Nb$ onto $x$.
\end{lemma}

The proof of this lemma uses the technique introduced in the construction
of $+$ in the proof of Proposition \ref{prop2.8}. Namely, we define the
graph of an $\iota$-function $g$ from $\Nb$ onto $x$ using $f$, $x$, and
$\Nb$ as parameters by writing a formula $\phi(z, f, x, \Nb)$ that asserts
``$z \in \Nb \times x$, and if $z = \langle a,b\rangle$ then there exists a
function $h: a+1 \to \Nb$ in $J_1$ such that $f(h(a)) = b$; if $k \leq h(a)$
satisfies $f(k) \in x$ then there exists $n \leq a$ such that
$f(h(n)) = f(k)$; and for each $n \leq a$ we have $f(h(n)) \in x$ and
$f(h(n)) \neq f(k)$ for any $k < h(n)$.'' The key point is that
$J_1$ contains the graphs of
all functions from $a+1 = \{0, \ldots, a\}$ into $\Nb$,
so that the required partial enumeration of $g$ is guaranteed to exist.

\begin{theo}\label{theorem2.10}
Let $x$ be an infinite $\iota$-set. Then there is an $\iota$-bijection
from $\Nb$ onto $x$.
\end{theo}

I omit details of the proof. The result can be established using the sets
$S_\alpha$ defined on p.\ 252 of \cite{Dev} with $\alpha = \omega + n$; one
can show inductively that (1) for every $n$ there is an $\iota$-surjection
from $\Nb$ onto $S_{\omega + n}$ and (2) for every $n$ there
exists $m$ such that $S_{\omega + m}$ contains the transitive
closure of $S_{\omega + n}$. Since $J_2 = \bigcup S_{\omega + n}$,
this plus Lemma \ref{lemma2.9} implies the theorem.

\begin{coro}\label{cor2.11}
(a) If $X$ is an $\iota$-class then so is
$$\Pc_{fin}(X) = \{y \in J_2: y \subseteq X\hbox{ is finite}\}.$$
If $x$ is an $\iota$-set then so is $\Pc_{fin}(x)$.

(b) If $X$ is an $\iota$-class then so is the set of all finite sequences
in $X$ (i.e., the graphs of all functions from $\{0, \ldots, n\}$ into $X$
for arbitrary $n \in \Nb$). If $x$ is an $\iota$-set then so is the set of
all finite sequences in $x$.
\end{coro}

\begin{proof}
(a) It is straightforward to verify that $\Pc_{fin}(X)$ is an $\iota$-class
for any $\iota$-class $X$. Now suppose $x$ is a set. The claim follows from
Proposition \ref{prop2.2} (d) if $x$ is finite, so suppose $x$ is infinite
and by Theorem \ref{theorem2.10} let $f: \Nb \to x$ be an $\iota$-bijection.
Observe that $\Pc_{fin}(\Nb) = \{x \in J_1: x \subseteq \Nb\}$ is an
$\iota$-set by the FDP. Then $\Pc_{fin}(x) = F_5(\Gamma(f), \Pc_{fin}(\Nb))$
where $F_5$ is the rudimentary function from \S \ref{sect1.1}, so
$\Pc_{fin}(x)$ is an $\iota$-set.

(b) Again, the statement about $\iota$-classes is straightforward.
To show that the set of finite sequences in $x$ is an $\iota$-set,
observe that $\Pc_{fin}(\Nb \times x)$ is an $\iota$-set by part (a);
as this contains the set of all finite sequences in $x$, the desired
result follows by the FDP.
\end{proof}

\subsection{Well-ordering and quotients}\label{sect2.4}

The following result is basic but its proof is a little involved.

\begin{theo}\label{theorem2.12}
There is an $\iota$-relation $\preceq_\Uc$ which well-orders $J_2$
in such a way that each initial segment is an $\iota$-set and each
$\iota$-set is contained in an initial segment.
\end{theo}

This result can be extracted from Lemma VI.2.7 of \cite{Dev} using
the fact mentioned above that the transitive closure of each
$S_{\omega + n}$ is contained in some $S_{\omega + m}$. I will call
$\preceq_\Uc$ the {\it universal well-ordering} on $J_2$.

Theorem \ref{theorem2.12} can be seen as a strong form of the axiom of choice
which yields a well-ordering not only of every $\iota$-set but of the
universal $\iota$-class. One basic way in which this is useful is in allowing
us to form quotients by equivalence relations. If we are given an equivalence
relation on an $\iota$-set which is an $\iota$-relation, we can form a
quotient $\iota$-set in the traditional way as the set of blocks of the
equivalence relation. But if we have an equivalence relation on a proper
$\iota$-class some of whose blocks are proper $\iota$-classes, then we clearly
cannot do this because proper $\iota$-classes cannot be elements of
$\iota$-classes. Instead, we can use the universal well-ordering to define
a version of the quotient by extracting a distinguished element from each
block.

\begin{defi}\label{def2.13}
An {\it equivalence $\iota$-relation} is an $\iota$-relation which is an
equivalence relation. Let $X$ be an $\iota$-class and let $\sim$ be an
equivalence $\iota$-relation on $X$. We define the {\it quotient} $X/{\sim}$
to be
\begin{eqnarray*}
X/{\sim} &=& \{a \in X:
(\forall b)[(b \in X \wedge a \sim b) \Rightarrow a \preceq_\Uc b]\}\\
&=& \{a \in J_2: \phi^{J_2}(a) \wedge
(\forall b \in J_2)[(\phi^{J_2}(b) \wedge
\psi^{J_2}(\langle a,b\rangle))
\Rightarrow \sigma^{J_2}(\langle a,b\rangle)]\}
\end{eqnarray*}
where $\phi$ is a formula that defines $X$, $\psi$ is a formula that
defines $\sim$, and $\sigma$ is a formula that defines $\preceq_\Uc$
(suppressing parameters).
\end{defi}

\begin{prop}\label{prop2.14}
The quotient of an $\iota$-set by an equivalence $\iota$-relation is an
$\iota$-set. The quotient of an $\iota$-class by an equivalence
$\iota$-relation is an $\iota$-class.
\end{prop}

\begin{proof}
Since $\preceq_\Uc$ is an $\iota$-relation its restriction to any $\iota$-set
is a small $\iota$-relation (Proposition \ref{prop2.5} (a), (c)). Thus, if $X$
and $\sim$ are $\iota$-sets then the definition of the quotient can be carried
out within $J_2$ and hence it is an $\iota$-set by the FDP.

Now let $X$ be an $\iota$-class and let $\sim$ be an equivalence
$\iota$-relation on $X$. Definability of $X/{\sim}$ was exhibited in
Definition \ref{def2.13}. For any $\iota$-set $x$, by Theorem
\ref{theorem2.12} we can find an $\iota$-set $y$ which contains $x$
and is an initial segment of the universal well-ordering. Then
$$(X/{\sim})\cap x = (X/{\sim}) \cap y \cap x
= [(X\cap y)/{\sim'}] \cap x$$
where $\sim'$ is the restriction of $\sim$ to $X \cap y$. Now $X \cap y$
is an $\iota$-set and $\sim'$ is an equivalence $\iota$-relation, so
$(X\cap y)/{\sim'}$ is an $\iota$-set by the first part of the proposition.
Its intersection with $x$ is then an $\iota$-set by Proposition \ref{prop2.2}
(b), and this shows that $X/{\sim}$ is an $\iota$-class.
\end{proof}

\section{The real line}\label{sect3}

\begin{quote}
{\it We now define the $\iota$-real line $\Rb_\iota$ in $J_2$. It is a
proper $\iota$-class, the standard ordering is an $\iota$-relation, and the
standard operations are $\iota$-functions. Various definitions of $\Rb_\iota$
are equivalent. In fact all of the standard classical functions from
$\Rb_\iota$ to $\Rb_\iota$ are $\iota$-functions.

We define $\iota$-open and $\iota$-closed $\iota$-classes in $\Rb_\iota$ and
establish their basic properties. The important notion of an $\iota$-set
being a ``proxy'' for an $\iota$-class is introduced here. We then define
$\iota$-compact and $\iota$-connected $\iota$-classes and discuss
$\iota$-continuous $\iota$-functions.}
\end{quote}

\subsection{Definition of $\Rb_\iota$}\label{sect3.1}
From here on I will usually invoke the FDP and the SDP
without giving any details; the reader should convince himself that these
uses are legitimate, which he should be able to do by verifying legitimacy
in detail in a few selected examples. In future sections I will use these
principles without even mentioning them.

The usual construction of $\Zb$ as the set of ordered pairs of natural
numbers modulo the equivalence $\langle m,n\rangle \sim \langle m',n'\rangle
\Leftrightarrow m + n' = m' + n$ can be straightforwardly carried
out in $J_2$, as can the usual construction of $\Qb$ in terms of
ordered pairs of integers modulo equivalence. There is likewise no obstacle
to defining the usual algebraic operations and order relation on $\Qb$.

\begin{defi}\label{def3.1}
A (lower) {\it Dedekind $\iota$-cut} is an $\iota$-subset of $\Qb$ which
is neither $\emptyset$ nor $\Qb$, has no greatest element, and contains
every element less than any element it contains. We also call a Dedekind
$\iota$-cut an {\it $\iota$-real}. The {\it $\iota$-real
line} $\Rb_\iota$ is the set of all Dedekind $\iota$-cuts. The functions
$+,\cdot: \Rb_\iota^2 \to \Rb_\iota$ and the ordering ${\leq} \subseteq
\Rb_\iota^2$ are defined in the standard way.
\end{defi}

We have the following basic facts.

\begin{theo}\label{theorem3.2}
$\Rb_\iota$ is an $\iota$-class, $+$ and $\cdot$ are $\iota$-functions
from $\Rb_\iota^2$ to $\Rb_\iota$, and $\leq$ is an $\iota$-relation on
$\Rb_\iota$. The subfield of $\Rb_\iota$ generated by any $\iota$-subset
is an $\iota$-set. Every nonempty $\iota$-subset of $\Rb_\iota$ that is
bounded above has a least upper bound and the map that takes each nonempty
$\iota$-subset that is bounded above to its least upper bound is an
$\iota$-function.
\end{theo}

\begin{proof}
$\Rb_\iota$ is an $\iota$-class by the SDP. Likewise $\leq_{\Rb_\iota}$ is
an $\iota$-relation by the SDP and $+$ and $\cdot$ are $\iota$-functions by
the SDP (used to check that they are $\iota$-classes) and the FDP (used to
check condition (2) on $\iota$-functions).

Given any nonempty $\iota$-subset $x \subseteq \Rb_\iota$, use Theorem
\ref{theorem2.10} to find an $\iota$-surjection $f: \Nb \to x$ and let $a$
be the $\iota$-set of (some nice encoding of) all words in the variables
$v_0, v_1, \ldots$ and the field operations. We can then use the FDP to show
that the set $b$ of ordered pairs $\langle r,q\rangle \in a \times \Qb$ such
that $q < r(f(0), f(1), \ldots)$ (i.e., $r$ evaluated with $v_n = f(n)$)
is an $\iota$-set. Finally, we can use the
rudimentary function $F(a) = \{G(y): y \in a\}$ from Proposition \ref{prop1.2}
with $G(y) = \{y\}$ to show that $\{ \{r\}: r \in a\}$ is an $\iota$-set;
applying the rudimentary function $F_5$ from \S \ref{sect1.1} to these two
$\iota$-sets yields that the subfield generated by $x$ is an $\iota$-set.

Next, $\bigcup_{y \in x} y$ is an $\iota$-set by Proposition \ref{prop2.2}
(b), and it belongs to $\Rb_\iota$ provided $x$ is bounded above. Since
$\leq_{\Rb_\iota}$ is just the inclusion relation, it follows that
$\bigcup_{y \in x} y$ is a least upper bound for $x$. Moreover, according
to Proposition \ref{prop1.2} the function $x \mapsto \bigcup_{y \in x} y$
is rudimentary and hence it is an $\iota$-function
by Proposition \ref{prop2.7}, so its restriction (call this $F$) to
$\Pc_\iota(\Rb_\iota)$ is an $\iota$-function by Proposition \ref{prop2.2} (d)
and Proposition \ref{prop2.6} (b); finally, the desired function
$F|_{F^{-1}[\Rb_\iota]}$ is an $\iota$-function by Proposition
\ref{prop2.6} (b) and (c).
\end{proof}

The other standard constructions of $\Rb_\iota$ are also available. For
example, a construction via Cauchy $\iota$-sequences (see Definition
\ref{def3.4}) can be carried out using
Definition \ref{def2.13}. One can prove the analog of Theorem \ref{theorem3.2}
here too. In fact, all standard constructions are equivalent by the following
result. Let an {\it $\iota$-field} be an $\iota$-class equipped with
$\iota$-functions which make it a field, such that the subfield generated
by any $\iota$-subset is an $\iota$-set. It is {\it $\iota$-ordered} if
it is an ordered field such that the partial order is an $\iota$-relation, and
it is {\it $\iota$-complete} if every nonempty $\iota$-subset that is bounded
above has a least upper bound and the map which takes each nonempty
$\iota$-subset that is bounded above to its least upper bound is an
$\iota$-function.

\begin{theo}\label{theorem3.3}
Every $\iota$-complete $\iota$-ordered $\iota$-field is isomorphic to
$\Rb_\iota$ via an $\iota$-bijection.
\end{theo}

The first part of the proof involves showing that every element of $\Fc$
corresponds to a Dedekind $\iota$-cut. To see this, let $x \in \Fc$. Then
the fact that $\leq_\Fc$ is an $\iota$-relation implies that $\{p \in \Qb_\Fc:
p <_\Fc x\}$ is an $\iota$-set, where $\Qb_\Fc$ is the canonical copy of $\Qb$
in $\Fc$. Since $\Qb_\Fc$ is the subfield of $\Fc$ generated by $1_\Fc$, it
is an $\iota$-set. Moreover, $\Qb_\Fc$ is isomorphic to $\Qb$ via an
$\iota$-bijection, so we infer that every element of $\Fc$ determines an
$\iota$-cut. Conversely, every $\iota$-cut determines an element of $\Fc$
since $\Fc$ is $\iota$-complete.
Density of $\Qb_\Fc$ in $\Fc$ is proven just as in the
classical case. This shows that there is an order-isomorphism between $\Fc$
and $\Rb_\iota$. The fact that this map respects $+$ and $\cdot$ is proven
just as in the classical case. Finally, the fact that the graph of the
isomorphism is an $\iota$-class can be shown using the FDP and the fact that
it is an $\iota$-function in both directions uses the hypothesis that the
map which takes each nonempty upper bounded $\iota$-set to its least upper
bound is an $\iota$-function.

The next result can be used to show that all of the standard functions from
$\Rb_\iota$ to $\Rb_\iota$ are $\iota$-functions. We use the following
terminology.

\begin{defi}\label{def3.4}
Let $X$ and $Y$ be $\iota$-classes. An {\it $\iota$-sequence} in $X$ is
an $\iota$-function from $\Nb$ to $X$. An {\it $\iota$-sequence of
$\iota$-functions} from $X$ to $Y$ is an $\iota$-function from
$\Nb \times X$ to $Y$.
\end{defi}

\begin{lemma}\label{lemma3.5}
If $(a_n)$ is a Cauchy $\iota$-sequence in $\Rb_\iota$ then it converges
to a limit in $\Rb_\iota$. The map that takes (graphs of) Cauchy
$\iota$-sequences in $\Rb_\iota$ to their limits is an $\iota$-function.
\end{lemma}

\begin{proof}
Regarding elements of $\Rb_\iota$ as Dedekind $\iota$-cuts, we have
$$\lim a_n = \bigcup_{n \in \Nb} \bigcap_{k \geq n} (a_k - 1/n).$$
Using this expression one can write a formula for $\lim a_n$ (as a
subset of $\Qb$) in terms of the graph of $(a_n)$, and it follows from the
FDP that the limit is an $\iota$-real. It is then fairly routine to verify
that the set of Cauchy $\iota$-sequences is an $\iota$-class and that the
map taking such a sequence to its limit is an $\iota$-function.
\end{proof}

\begin{prop}\label{prop3.6}
Let $X \subseteq \Rb_\iota$ be an $\iota$-class. Then any pointwise limit
of an $\iota$-sequence of $\iota$-functions from $X$ to $\Rb_\iota$ is an
$\iota$-function.
\end{prop}

\begin{proof}
Let $F: \Nb \times X \to \Rb_\iota$ be an $\iota$-sequence of
$\iota$-functions and suppose that
for each $a \in X$ the $\iota$-sequence $(F(\cdot, a))$ is Cauchy. We claim
that the map $G$ which takes $a \in X$ to the graph of the $\iota$-sequence
$(F(\cdot, a))$ is an $\iota$-function. Definability of its graph is
straightforward. To show that $\Gamma(G)$ is an $\iota$-class it is sufficient
to show that its intersection with $x \times J_2$ is an $\iota$-set for every
$\iota$-subset $x$ of $X$, since for any $\iota$-set $x$ we have
$$\Gamma(G) \cap x = \big(\Gamma(G) \cap ({\rm dom}(x) \times J_2)\big)
\cap x.$$
So given an $\iota$-set $x \subseteq X$, find an $\iota$-subset $y$ of
$\Rb_\iota$ such that $F[\Nb \times x] \subseteq y$; then
$\Gamma(G) \cap (x \times J_2)$ is an $\iota$-subset of $x \times J_2$ by
the FDP, using the $\iota$-set $\Gamma(F) \cap (\Nb \times x \times y)$ as
a parameter. The second condition on $\iota$-functions is proven similarly
and we conclude that $G$ is an $\iota$-function. The desired result then
follows from the lemma since any composition of $\iota$-functions is an
$\iota$-function (Proposition \ref{prop2.6} (d)).
\end{proof}

\begin{coro}\label{cor3.7}
All of the standard functions on $\Rb_\iota$ are $\iota$-functions.
\end{coro}

All polynomials on $\Rb_\iota$ are $\iota$-functions since $+$ and $\cdot$
are $\iota$-functions (Theorem \ref{theorem3.2}), the constant functions and
the identity map on $\Rb_\iota$ are $\iota$-functions (Proposition
\ref{prop2.6} (e)), and compositions of $\iota$-functions are
$\iota$-functions (Proposition \ref{prop2.6} (d)). All standard continuous
functions ($\sin t$, $\cos t$, $\ln t$, $e^t$, etc.) are $\iota$-functions
by Proposition \ref{prop3.6} together with this observation. Standard
discontinuous functions such as the jump function $f(t) = (0$ for $t \leq 0$
and $1$ for $t > 0)$ can also be seen to be $\iota$-functions via pointwise
approximation by polynomials, or more simply by a direct proof.

\subsection{Open $\iota$-classes}\label{sect3.2}

We now consider $\iota$-classes $U \subseteq \Rb_\iota$ which are ``open''. No
such $\iota$-class will be an $\iota$-set unless it is empty, because every
nonempty open $\iota$-class will contain an interval of positive length, and
no such interval can be contained in an $\iota$-set. The right definition seems
to be the following. Let $\Qb^+ = \Qb \cap (0, \infty)$ and $\Rb_\iota^+ =
\Rb_\iota \cap (0,\infty)$.

\begin{defi}\label{def3.8}
An $\iota$-class $U \subseteq \Rb_\iota$ is {\it $\iota$-open} if there is
an $\iota$-set $u \subseteq \Qb \times \Qb^+$ such that for any $\iota$-real
number $r \in \Rb_\iota$ we have $r \in U$ if and only if there exists
$\langle p,q\rangle \in u$ with $|r-p| < q$. We call $u$ a {\it proxy}
for $U$.
\end{defi}

That is, every $\iota$-open $\iota$-class is a union of an $\iota$-set of
balls with rational centers and rational radii. Note that this is {\it not}
the same as saying that for every $r \in U$ there exists such a ball containing
$r$ and contained in $U$. The latter condition is weaker since we require
that $u$ be an $\iota$-set (but see \S \ref{sect3.4}).

The concept of a proxy is useful if we want to make a statement about an
$\iota$-set of $\iota$-open $\iota$-classes. This does not make literal
sense because a proper $\iota$-class can never be a member of another
$\iota$-set or $\iota$-class, but we can often give such a statement a
reasonable meaning using proxies. For example, in part (a) of the next
result we assert that the union of any $\iota$-set of $\iota$-open
$\iota$-classes is $\iota$-open; this should be understood as abbreviating
the assertion ``given an $\iota$-set of proxies of $\iota$-open
$\iota$-classes, the union of the corresponding $\iota$-open $\iota$-classes
is $\iota$-open''.

\begin{prop}\label{prop3.9}

(a) The union of any $\iota$-set of $\iota$-open $\iota$-classes is an
$\iota$-open $\iota$-class.

(b) Any intersection of finitely many $\iota$-open $\iota$-classes is an
$\iota$-open $\iota$-class.

(c) Let $u \subseteq \Qb \times \Qb^+$ be an $\iota$-set. Then
$$U = \{r \in \Rb_\iota: |r-p| < q\hbox{ for some }\langle p,q\rangle \in u\}$$
is an $\iota$-open $\iota$-class.
\end{prop}

\begin{proof}

(a) This holds because if $x$ is an $\iota$-set of proxies of $\iota$-open
$\iota$-classes then $\bigcup_{y \in x} y$ is a proxy for the union of the
corresponding $\iota$-open $\iota$-classes.

(b) The intersection of any two balls with rational radii and rational
centers is again a ball with a rational radius and a rational center.
Thus if $U$ and $V$ are two $\iota$-open $\iota$-classes, a proxy for their
intersection is obtained by intersecting every ball in a proxy for $U$
with every ball in a proxy for $V$. This is enough.

(c) $U$ is an $\iota$-class by the SDP, and it is then $\iota$-open by
definition.
\end{proof}

\begin{coro}\label{cor3.10}
Let $u \subseteq \Rb_\iota \times \Rb_\iota^+$ be an $\iota$-set. Then
$$U = \{r \in \Rb_\iota: |r-p| < q\hbox{ for some }\langle p,q\rangle \in u\}$$
is an $\iota$-open $\iota$-class.
\end{coro}

\begin{proof}
Define
$$v = \{\langle p,q\rangle \in \Qb \times \Qb^+:
q + |p-r| \leq s\hbox{ for some }\langle r,s \rangle \in u\}.$$
This is an $\iota$-set by the FDP. By Proposition \ref{prop3.9} (c),
$v$ is a proxy for an $\iota$-open $\iota$-class, which by inspection
equals $U$.
\end{proof}

In particular, the {\it $\iota$-open ball}
$$B_q(p) = \{r \in \Rb_\iota: |r - p| < q\}$$
is an $\iota$-open $\iota$-class for any $p \in \Rb_\iota$ and
$q \in \Rb_\iota^+$.

\subsection{Closed $\iota$-classes}\label{sect3.3}

\begin{defi}\label{def3.11}
Let $c$ be an $\iota$-subset of $\Rb_\iota$. The {\it $\iota$-closure} of
$c$ is the set of all $\iota$-real numbers that are limits of Cauchy
$\iota$-sequences in $c$. An $\iota$-class $C \subseteq \Rb_\iota$ is
{\it $\iota$-closed} if it is the $\iota$-closure of an $\iota$-set $c$,
in which case $c$ is a {\it proxy} for $C$.
\end{defi}

\begin{prop}\label{prop3.12}
Any Cauchy $\iota$-sequence in an $\iota$-closed $\iota$-class converges
to an element of that $\iota$-class.
\end{prop}

\begin{proof}
Let $C$ be an $\iota$-closed $\iota$-class with proxy $c$ and let $(a_n)$ be
a Cauchy $\iota$-sequence in $C$. We can define a Cauchy sequence $(b_n)$ with
values in $c$ that has the same limit by letting $b_n$ be the first element of
$c$, with respect to $\preceq_\Uc$, whose distance to $a_n$ is less than $1/n$.
This is an $\iota$-sequence because its graph is an $\iota$-subset of
$\Nb \times c$ by the FDP. We conclude that $\lim a_n$ belongs to $C$.
\end{proof}

\begin{theo}\label{theorem3.13}
An $\iota$-class $C \subseteq \Rb_\iota$ is $\iota$-closed if and only if
$\Rb_\iota - C$ is $\iota$-open.
\end{theo}

\begin{proof}
Suppose $C \subseteq \Rb_\iota$ is $\iota$-closed and let $c$ be a proxy
for $C$. Then
$$u = \{\langle p,q\rangle \in \Qb \times \Qb^+:
q \leq |p-r|\hbox{ for all }r \in c\}$$
is an $\iota$-set by the FDP. Let $U$ be the corresponding $\iota$-open
$\iota$-class. It is clear that $U \subseteq \Rb_\iota - C$. Conversely,
for any $r \in \Rb_\iota - C$ there must exist $\epsilon > 0$ such that
$B_\epsilon(r) \cap c = \emptyset$; otherwise we could find a Cauchy
$\iota$-sequence in $c$ converging to $r$ by the method used to construct
$(b_n)$ in Proposition \ref{prop3.12}. Density of $\Qb$ in $\Rb_\iota$ now
implies $r \in U$. We conclude that $\Rb_\iota - C = U$ is $\iota$-open.

Now let $U \subseteq \Rb_\iota$ be $\iota$-open; we must show that its
complement is $\iota$-closed. We may assume $U$ is nonempty.
Let $u$ be a proxy for $U$. Then let $x$
be the set of ordered pairs of rationals $\langle p,q\rangle$ with the
property that $p < q$ and the interval $[p,q]$ is covered by finitely
many balls indexed by $u$. That is, $\langle p,q \rangle \in x$ if and
only if $p < q$ and
there exists a finite sequence $\langle p_1, q_1\rangle, \ldots,
\langle p_n, q_n\rangle$ in $u$ such that (1) $|p - p_1| < q_1$;
(2) $|p_i - p_{i+1}| < q_i + q_{i+1}$ for $1 \leq i < n$; and
(3) $|p_n - q| < q_n$. By the FDP, $x$ is an $\iota$-set. We can also
use the FDP to show that
$$y = \{p \in \Qb: |p' - p| < q'\hbox{ for some }\langle p',q'\rangle \in u\}$$
is an $\iota$-set. Next, let $x'$ be the $\iota$-set
$$x' = x \cup \{\langle p,q\rangle \in \Qb^2: p \in y\hbox{ and }q \leq p\}
\cup \{\langle p,q\rangle \in \Qb^2: p \not\in y\hbox{ and }q < p\}.$$
Some simple manipulations using rudimentary functions (in particular the
function $F_5$ from \S \ref{sect1.1}) then yield an $\iota$-set $c$
whose elements are the $\iota$-cuts $\{q \in \Qb: \langle p,q \rangle
\in x'\}$ (for all $p \in \Qb$ such that this set is properly contained
in $\Qb$). We have constructed $c$ so that it contains all rationals in
the complement of $U$, as well as, for each rational $p \in U$, the
smallest irrational not in $U$ that is greater than $p$. It is now a
routine matter to use the technique of the proof of Proposition
\ref{prop3.12} to check that the closure of $c$ equals $\Rb_\iota - U$
(see Proposition \ref{prop3.15} below).
\end{proof}

\begin{coro}\label{cor3.14}

(a) The intersection of any $\iota$-set of $\iota$-closed $\iota$-classes
is an $\iota$-closed $\iota$-class.

(b) Any union of finitely many $\iota$-closed $\iota$-classes is an
$\iota$-closed $\iota$-class.

(c) The $\iota$-closure of any $\iota$-subset of $\Rb_\iota$ is an
$\iota$-closed $\iota$-class.
\end{coro}

\begin{proof}
(a) Given an $\iota$-set of proxies for $\iota$-closed $\iota$-classes,
the construction of $U$ from $c$ in the proof of Theorem \ref{theorem3.13}
yields an $\iota$-set of proxies for the complementary $\iota$-open
$\iota$-classes. (Assuming $x$ is an infinite $\iota$-set of proxies for
$\iota$-closed $\iota$-classes, use an $\iota$-bijection between $\Nb$ and
$x$ to find an $\iota$-subset of $\Nb \times \Qb \times \Qb^+$ whose
restriction to each $n \in \Nb$ is a proxy for the $\iota$-open $\iota$-class
that is complementary to the corresponding $\iota$-closed $\iota$-class in $x$.
Then use the rudimentary function $F_5$ from \S \ref{sect1.1} to get the
desired $\iota$-set of proxies for the complementary $\iota$-open
$\iota$-classes.) The result now follows from Proposition \ref{prop3.9}
(a) and Theorem \ref{theorem3.13}.

(b) This follows directly from Proposition \ref{prop3.9} (b) and
Theorem \ref{theorem3.13}.

(c) Let $c \subseteq \Rb_\iota$ be an $\iota$-set. Then the proof of Theorem
\ref{theorem3.13} shows that the complement of its $\iota$-closure is the
union of an $\iota$-set of balls with rational centers and radii.
The latter is an $\iota$-open $\iota$-class by Proposition \ref{prop3.9} (c),
so the $\iota$-closure of $c$ is an $\iota$-closed $\iota$-class by Theorem
\ref{theorem3.13}.
\end{proof}

\subsection{Separable $\iota$-classes}\label{sect3.4}

As I mentioned at the beginning of \S \ref{sect3.2}, for an $\iota$-subclass
of $\Rb_\iota$ to be $\iota$-open it is not sufficient that every element of
the $\iota$-class be contained in a ball contained in the $\iota$-class. For
example, it is possible to construct a proper $\iota$-class that contains
exactly one $\iota$-real number in the interval $[n,n+1]$ for each $n \in \Nb$.
Although it is closed under limits of Cauchy $\iota$-sequences (since every
Cauchy $\iota$-sequence is eventually constant), this will not be an
$\iota$-closed $\iota$-class according to Definition \ref{def3.11} because
it is not the $\iota$-closure of an $\iota$-set. Its complement will be an
$\iota$-class that is not $\iota$-open according to Definition \ref{def3.8}
but which does have the ball property mentioned above.

However, the familiar equivalences do hold if we make additional separability
assumptions.

\begin{prop}\label{prop3.15}
Let $X \subseteq \Rb_\iota$ be an $\iota$-class and let $x \subseteq X$ be
an $\iota$-set. Then the following are equivalent:

(a) Every $\iota$-open ball about a point of $X$ intersects $x$.

(b) Every $\iota$-open ball with rational center and radius that intersects $X$
also intersects $x$.

(c) Every $\iota$-open $\iota$-class that intersects $X$ also intersects $x$.

(d) $X$ is contained in the $\iota$-closure of $x$.
\end{prop}

\begin{proof}
The equivalence of (a), (b), and (c) is trivial, as is the fact that
(d) implies them. The final implication is proven by constructing a
Cauchy $\iota$-sequence in $x$ that converges to a given point in $X$
by the technique used in the proof of Proposition \ref{prop3.12}.
\end{proof}

\begin{defi}\label{def3.16}
Let $X \subseteq \Rb_\iota$ be an $\iota$-class and let $x \subseteq X$ be
an $\iota$-set. We say that $x$ is {\it $\iota$-dense} in $X$ if any of the
equivalent conditions of Proposition \ref{prop3.15} holds. We say that $X$
is {\it $\iota$-separable} if it contains an $\iota$-dense $\iota$-subset.
\end{defi}

\begin{lemma}\label{lemma3.17}
Every $\iota$-open $\iota$-class is $\iota$-separable, as is every
$\iota$-closed $\iota$-class.
\end{lemma}

\begin{prop}\label{prop3.18}
Let $X \subseteq \Rb_\iota$ be an $\iota$-separable $\iota$-class and
let $Y = \Rb_\iota - X$.

(a) $X$ is $\iota$-closed if and only if it is closed under limits of all
Cauchy $\iota$-sequences in $X$.

(b) $Y$ is $\iota$-open if and only if for every $r \in Y$ there exists
$\epsilon > 0$ such that $B_\epsilon(r) \subseteq Y$.
\end{prop}

\begin{proof}
(a) The forward direction was Proposition \ref{prop3.12}. For the reverse
direction, let $c$ be an $\iota$-dense $\iota$-subset of $X$; then $X$ is
automatically contained in the $\iota$-closure of $c$, and it contains the
$\iota$-closure of $c$ because it is closed under limits of Cauchy
$\iota$-sequences. So $X$ is $\iota$-closed.

(b) The forward direction is easy. For the reverse direction, let $c$ be
an $\iota$-dense $\iota$-subset of $X$ and suppose every element of $Y$ is
contained in a ball contained in $Y$. Let $C$ be the $\iota$-closure of $c$;
it is then easy to verify that $Y = \Rb_\iota - C$, which implies that $Y$
is $\iota$-open by Theorem \ref{theorem3.13}.
\end{proof}

\subsection{Compactness and connectedness}\label{sect3.5}

\begin{defi}\label{def3.19}
Let $K \subseteq \Rb_\iota$ be an $\iota$-class. We say that
$K$ is {\it $\iota$-compact} if any $\iota$-set of $\iota$-open $\iota$-classes
which covers $K$ has a finite subcover. It is {\it $\iota$-connected} if
there do not exist $\iota$-open $\iota$-classes $U$ and $V$ such that
$K \subseteq U \cup V$, $K \cap U \neq \emptyset$, $K \cap V \neq \emptyset$,
and $U \cap V = \emptyset$.
\end{defi}

Recall our convention about phrases like ``$\iota$-set of $\iota$-open
$\iota$-classes'' from \S \ref{sect3.3}: the above really means that
whenever $x$ is an $\iota$-set of proxies of $\iota$-open $\iota$-classes
and the union of the corresponding $\iota$-classes contains $K$, there is
a finite subset of $x$ with the same property.

\begin{theo}\label{theorem3.20}
Let $K$ be an $\iota$-separable $\iota$-subclass of $\Rb_\iota$. Then the
following are equivalent:

(i) $K$ is an $\iota$-closed and bounded $\iota$-class;

(ii) $K$ is $\iota$-compact;

(iii) $K$ is bounded and contains the limits of all of its Cauchy
$\iota$-sequences;

(iv) every $\iota$-sequence in $K$ has an $\iota$-subsequence which
converges to a limit in $K$.
\end{theo}

\begin{proof}

(i) $\Rightarrow$ (ii): This reduces to the assertion that $[0,1]$ is
$\iota$-compact in the usual way. Since every $\iota$-open $\iota$-class
is a union of an $\iota$-set of balls with rational centers and radii,
we further reduce to covers of $[0,1]$ by such balls. The proof of
$\iota$-compactness is now the standard one, already essentially
given in the last part of the proof of Theorem \ref{theorem3.13}.

(ii) $\Rightarrow$ (iii): Suppose $K$ is $\iota$-compact. If $K$ were not
bounded then $\{B_n(0): n \in \Nb\}$ would be an $\iota$-set of $\iota$-open
$\iota$-classes which covers $K$ but has no finite subcover. So $K$ is
bounded. Similarly, if there were a Cauchy $\iota$-sequence $(r_n)$ in $K$
whose limit did not belong to $K$ then $\{B_q(p): p,q \in \Qb$ and $q <
|p-r|/2\}$ would contradict $\iota$-compactness,
where $r = \lim r_n$. So $K$ is closed under convergence of Cauchy
$\iota$-sequences.

(iii) $\Rightarrow$ (iv): Suppose $K$ is bounded and contains the limits of
all of its Cauchy $\iota$-sequences and let $(r_n)$ be an $\iota$-sequence
in $K$. For each $k$ let $I_k$ be the leftmost interval of the form
$[j/2^k, (j+1)/2^k]$ which contains infinitely many terms of the
$\iota$-sequence $(r_n)$ and let $n_k$ be the first index greater than
$n_{k-1}$ such that $r_{n_k} \in I_k$. Then $(r_{n_k})$ is a Cauchy
$\iota$-subsequence of $(r_n)$ since the $I_k$ are nested, and hence
it converges to a limit in $K$.

(iv) $\Rightarrow$ (i) Suppose every
$\iota$-sequence in $K$ has an $\iota$-subsequence which converges to a
limit in $K$ and let $c$ be an $\iota$-dense $\iota$-subset of $K$. If $c$
were unbounded we could find an $\iota$-sequence in $c$ which diverges to
$\pm \infty$ and hence has no convergent $\iota$-subsequence, so $c$ must
be bounded. This implies that $K$ is bounded. Now $K$ is automatically
contained in the $\iota$-closure of $c$; conversely, any Cauchy
$\iota$-sequence in $c$ must converge to a limit in $K$ because by
assumption it has an $\iota$-subsequence which converges to a limit in $K$.
Thus $K$ is the $\iota$-closure of $c$ and hence $K$ is $\iota$-closed.
\end{proof}

(Note that only the proof of (iv) $\Rightarrow$ (i) used $\iota$-separability.)

The characterization of $\iota$-connected $\iota$-classes in $\Rb_\iota$
is proven in just the same way as in the classical case.

\begin{theo}\label{theorem3.21}
An $\iota$-subclass of $\Rb_\iota$ is $\iota$-connected if and only if
it is an interval.
\end{theo}

\subsection{Continuous $\iota$-functions}\label{sect3.6}

\begin{defi}\label{def3.22}
Let $X$ be an $\iota$-subclass of $\Rb_\iota$. An $\iota$-function
$F: X \to \Rb_\iota$ is {\it $\iota$-continuous} if the inverse image of
every $\iota$-open $\iota$-class in $\Rb_\iota$ is the intersection of an
$\iota$-open $\iota$-class with $X$.
\end{defi}

\begin{theo}\label{theorem3.23}
Let $X$ be an $\iota$-separable $\iota$-subclass of $\Rb_\iota$ and let
$F: X \to \Rb_\iota$ be an $\iota$-function. Then the following are equivalent:

(a) $F$ is $\iota$-continuous.

(b) The inverse image of every $\iota$-closed $\iota$-class under $F$ is the
intersection of an $\iota$-closed $\iota$-class with $X$.

(c) For any $\iota$-set $c \subseteq X$ with $\iota$-closure $C$, the
$\iota$-closure of $F[c]$ contains $F[C \cap X]$.

(d) $F$ preserves convergence of $\iota$-sequences.

(e) For every $r \in X$ and every $\epsilon > 0$ there exists
$\delta > 0$ such that $|r - s| < \delta$ implies $|F(r) - F(s)| <
\epsilon$.

If any of these conditions holds then there is an $\iota$-function which takes
any proxy $u$ of an $\iota$-open $\iota$-class $U \subseteq \Rb_\iota$ to a
proxy $v$ of an $\iota$-open $\iota$-class $V \subseteq \Rb_\iota$ such that
$F^{-1}[U] = X \cap V$.
\end{theo}

\begin{proof}
(a) $\Rightarrow$ (b) $\Rightarrow$ (c) $\Rightarrow$ (d): The classical
proofs carry over without alteration.

(d) $\Rightarrow$ (e): Suppose (e) fails and choose
$r$, $\epsilon$ such that for every $\delta > 0$ there exists $s$
with $|r - s| < \delta$ and $|F(r) - F(s)| \geq \epsilon$. This implies
that there is a sequence in $X$ which converges to $r$ but whose image
does not converge to $F(r)$. However, it is not obvious that there is
an $\iota$-sequence with this property so the contradiction is not immediate.

Let $x$ be an $\iota$-dense $\iota$-subset of $X$. If for every $k \in \Nb$
there exists $s \in x$ such that $|r - s| < 1/k$ and $|F(r) - F(s)| \geq
\epsilon/2$ then we can find an $\iota$-sequence which converges to $r$
whose image does not converge to $F(r)$. This would falsify (d) and we
would be done. Therefore we now assume that there exists $k \in \Nb$ such that
$|F(r) - F(s)| < \epsilon/2$ for every $s \in x$ such that $|r - s| < 1/k$.

By the choice of $r$ we know there exists $s \in X$ with $|r - s| < 1/k$
and $|F(r) - F(s)| \geq \epsilon$. But we also know that $x$ is $\iota$-dense
in $X$ and hence there is an $\iota$-sequence $(s_n)$ in $x$ that converges
to $s$, and $|r - s_n| < 1/k$ for sufficiently large $n$. But then the
choice of $k$ yields $|F(r) - F(s_n)| < \epsilon/2$ for sufficiently large
$n$, which implies that $F(s_n) \not\to F(s)$, and we have falsified (d).

(e) $\Rightarrow$ (a): Suppose (e) holds and let $U \subseteq \Rb_\iota$ be an
$\iota$-open $\iota$-class with proxy $u$; we must show that $F^{-1}(U)$ is
the intersection of an $\iota$-open $\iota$-class with $X$. Let $x$ be an
$\iota$-dense $\iota$-subset of $X$ and let $v$ be the $\iota$-set of pairs
$\langle p,q\rangle \in \Qb\times\Qb^+$ such that $F(x \cap B_q(p))
\subseteq B_{q' - 1/k}(p')$ for some $\langle p',q'\rangle \in u$ and some
$k \in \Nb$. We claim that $v$ is a proxy for an $\iota$-open $\iota$-class
$V$ whose intersection with $X$ equals $F^{-1}(U)$.

First we show that $V$ contains $F^{-1}(U)$. To see this let $r \in
F^{-1}(U)$. Find $\langle p',q'\rangle \in u$ and $k \in \Nb$ such
that $F(r) \in B_{q'-1/k}(p')$ and let $\epsilon = q' - 1/k - |p' - F(r)|$.
By (e) we can then find
$\delta > 0$ such that $|r - s| < \delta$ implies $|F(r) - F(s)| < \epsilon$.
Finally we can find $\langle p,q\rangle$ such that
$r \in B_q(p) \subseteq B_\delta(r) \subseteq B_{q' - 1/k}(p')$, which
implies that $\langle p,q\rangle \in v$, and we conclude that $r \in V$.
This shows one containment.

Now we must show that $F^{-1}(U)$ contains $V \cap X$. Suppose not and
fix $r \in V \cap X$ such that $F(r) \not\in U$. Fix $\langle p,q\rangle
\in v$ such that $r \in B_q(p)$; by (e) and density of $x$ in $X$ we can
then find, for every $k \in \Nb$, an $s \in x \cap B_q(p)$ such that
$|F(r) - F(s)| < 1/k$. But since $F(r) \not\in U$ this
contradicts the fact that $F(x \cap B_q(p))$ is contained in some
$B_{q' - 1/k}(p')$. We conclude that $V \cap X \subseteq F^{-1}(U)$.
This completes the proof that (e) implies (a).

In the proof of (e) $\Rightarrow$ (a) the proxy $v$ is definable
from the proxy $u$ and this shows that the map $u \mapsto v$ is an
$\iota$-function by the FDP.
\end{proof}

The proof of the next result is no different from the classical case.

\begin{theo}\label{theorem3.24}
The composition, sum, and product of any two $\iota$-continuous
$\iota$-functions is $\iota$-continuous.
\end{theo}

\begin{theo}\label{theorem3.25}
(a) The image of an $\iota$-separable $\iota$-compact $\iota$-class under
an $\iota$-continuous $\iota$-function is an $\iota$-separable
$\iota$-compact $\iota$-class.

(b) The image of an $\iota$-connected $\iota$-class under an $\iota$-continuous
$\iota$-function is an $\iota$-connected $\iota$-class.
\end{theo}

\begin{proof}
(a) Let $K$ be an $\iota$-separable $\iota$-compact $\iota$-class and let $F$
be an $\iota$-continuous $\iota$-function. One proves the $\iota$-compactness
property of $F[K]$ just as in the classical case, using the last part of
Theorem \ref{theorem3.23} to pull an $\iota$-open cover of $F[K]$ back to
an $\iota$-open cover of $K$. However, we must also show that $F[K]$ is
an $\iota$-separable $\iota$-class. Let $c$ be an $\iota$-dense $\iota$-subset
of $K$; we claim that $F[K]$ equals the $\iota$-closure of $F[c]$. The forward
containment follows from Theorem \ref{theorem3.23} (c). For the
reverse containment, let $r$ belong to the $\iota$-closure of $F[c]$.
If $r \not\in F[K]$ then $\{B_q(p): p,q \in \Qb$ and $q < |p-r|/2\}$
would be an $\iota$-open cover of $F[K]$ with no finite
subcover, which would pull back to an $\iota$-open cover of $K$ with
no finite subcover. This contradicts $\iota$-compactness of $K$ and
we conclude that $F[K]$ is the closure of $F[c]$, so it is an
$\iota$-class by Corollary \ref{cor3.14} (c).

(b) Again, one proves $\iota$-connectedness of the image just as in
the classical case. This does not use the fact that the image is an
$\iota$-class, but it implies that the image is an interval and hence
it is in fact an $\iota$-class.
\end{proof}

\begin{coro}\label{cor3.26}
Any $\iota$-continuous $\iota$-function on an $\iota$-separable
$\iota$-compact $\iota$-subclass of $\Rb_\iota$ is bounded and
achieves its maximum and minimum. Any
$\iota$-continuous $\iota$-function on an $\iota$-connected $\iota$-subclass
of $\Rb_\iota$ achieves all intermediate values.
\end{coro}

\section{Topology}\label{sect4}

\begin{quote}
{\it We define $\iota$-metric and $\iota$-topological spaces and discuss
their basic properties.}
\end{quote}

\subsection{Metric spaces}\label{sect4.1}

Most of the following material on metric spaces is a straightforward
generalization of material in Section \ref{sect3}.

\begin{defi}\label{def4.1}
An {\it $\iota$-metric space} is an $\iota$-class $X$ together with an
$\iota$-function $D: X \times X \to [0,\infty)$ which satisfies the usual
metric axioms, such that the map that takes convergent $\iota$-sequences
to their limits is an $\iota$-function. It is {\it $\iota$-separable} if
it contains an $\iota$-subset $x$ which intersects every ball and it is
{\it $\iota$-complete} if every Cauchy $\iota$-sequence converges.
\end{defi}

\begin{defi}\label{def4.2}
Let $X$ be an $\iota$-metric space. We define the {\it completion} of $X$ to
be the $\iota$-class of Cauchy $\iota$-sequences in $X$ modulo the standard
equivalence, with distance $\iota$-function defined in the usual way.
\end{defi}

\begin{prop}\label{prop4.3}
The completion of any $\iota$-metric space $X$ is an $\iota$-complete
$\iota$-metric space. The canonical embedding of $X$ in its completion
is an $\iota$-bijection between $X$ and its image.
\end{prop}

In Definition \ref{def4.1} we included an extra assumption stating that
the map that takes convergent $\iota$-sequences to their limits is an
$\iota$-function. This assumption is harmless because the standard
construction of the $\iota$-completion ensures this condition. That is,
if we apply the $\iota$-completion construction to an $\iota$-class $X$
equipped with any $\iota$-function $D: X \times X \to [0,\infty)$ which
satisfies the usual metric axioms, the result is an $\iota$-complete
$\iota$-metric space. The image of $X$ in its completion will also be
an $\iota$-metric space, though if $X$ is not an $\iota$-metric space
then the canonical embedding will not be an $\iota$-bijection of $X$
with its image because the inverse map will not be an $\iota$-function.

\begin{prop}\label{prop4.4}
In any $\iota$-metric space, any pointwise limit of an $\iota$-sequence
of $\Rb_\iota$-valued $\iota$-functions is an $\iota$-function.
\end{prop}

\begin{defi}\label{def4.5}
Let $X$ be an $\iota$-metric space. An $\iota$-class $U \subseteq X$
is {\it $\iota$-open} if there is an $\iota$-set $u \subseteq X \times
\Rb_\iota^+$ such that for any $r \in X$ we have $r \in U$ if and only
if there exists $\langle p,q \rangle \in u$ with $D(r,p) < q$. An
$\iota$-class $C \subseteq X$ is {\it $\iota$-closed} if its complement
is $\iota$-open. We call $u$ a {\it proxy} both for the $\iota$-open
$\iota$-class $U$ and the $\iota$-closed $\iota$-class $X - U$.
\end{defi}

\begin{theo}\label{theorem4.6}
Let $X$ be an $\iota$-separable $\iota$-metric space.

(a) The union of any $\iota$-set of $\iota$-open $\iota$-classes is
$\iota$-open and the intersection of any $\iota$-set of $\iota$-closed
$\iota$-classes is $\iota$-closed.

(b) Any intersection of finitely many $\iota$-open $\iota$-classes is
$\iota$-open and any union of finitely many $\iota$-closed $\iota$-classes
is $\iota$-closed.

(c) Every $\iota$-subset of $X \times \Rb_\iota^+$ is a proxy for an
$\iota$-open class (and for its complementary $\iota$-closed class).
\end{theo}

(We need $\iota$-separability of $X$ in part (b), and also to ensure
that $X$ is $\iota$-open.)

\begin{defi}\label{def4.7}
The {\it $\iota$-closure} of an $\iota$-class $Y$ in an $\iota$-metric
space $X$ is the set of all limits of convergent $\iota$-sequences in $Y$.
$Y$ is {\it $\iota$-dense} in $X$ if its $\iota$-closure equals $X$.
\end{defi}

\begin{defi}\label{def4.8}
(a) An $\iota$-class $K$ contained in an $\iota$-metric space is
{\it $\iota$-totally bounded} if there is an $\iota$-function
$f: \Nb \to \Pc_{fin}(K)$ such that for each $k \in \Nb$ the
$\iota$-class $K$ is covered by the balls of radius $2^{-k}$ about the
elements of $f(k)$.

(b) An $\iota$-subclass of an $\iota$-metric space is
{\it $\iota$-compact} if any $\iota$-set of $\iota$-open $\iota$-classes
which covers it has a finite subcover.
\end{defi}

\begin{theo}\label{theorem4.9}
Let $X$ be an $\iota$-separable $\iota$-metric space. The following are
equivalent:

(i) $X$ is $\iota$-compact;

(ii) $X$ is $\iota$-complete and $\iota$-totally bounded;

(iii) Every $\iota$-sequence in $X$ has a convergent $\iota$-subsequence.
\end{theo}

\begin{proof}
(i) $\Rightarrow$ (ii): Suppose $X$ is not $\iota$-complete and let $(r_n)$
be a Cauchy $\iota$-sequence with no limit. Let $x$ be an $\iota$-dense
$\iota$-subset of $X$ and let $y \subseteq x \times \Qb^+$ be the $\iota$-set
of pairs $\langle p,q\rangle$ such that $2q < \lim D(p,r_n)$. Then it is easy
to see that $\{B_q(p): \langle p,q\rangle \in y\}$ is an $\iota$-set of
$\iota$-open balls that covers $X$ but has no finite subcover. We conclude
that any $\iota$-compact space is $\iota$-complete.

Now suppose $X$ is $\iota$-compact. By Theorem \ref{theorem2.10} we may write
$x = \{r_n: n \in \Nb\}$ for some $\iota$-sequence $(r_n)$. We claim that
for each $k \in \Nb$ there exists $m \in \Nb$ such that every $r \in x$
satisfies $D(r,r_n)<2^{-k-1}$ for some $n \leq m$. If such an $m$ exists for
each $k$ then we can verify total boundedness by setting $f(k) =
\{r_1, \ldots, r_m\}$. But if there were no such $m$ for some value of $k$
then the balls $B_{2^{-k-1}}(r_n)$ would cover $X$ yet have no finite
subcover, contradicting $\iota$-compactness. So $X$ is $\iota$-totally bounded.

(ii) $\Rightarrow$ (iii): Suppose $X$ is $\iota$-complete and $\iota$-totally
bounded and let $(r_n)$ be an $\iota$-sequence in $X$. Let $f$ be an
$\iota$-function which verifies total boundedness. For each $k$ let $g(k)$
be the set of $r \in f(k)$ such that $B_{2^{-k}}(r)$ contains infinitely many
terms of the $\iota$-sequence $(r_n)$. Observe that $g(k)$ is nonempty for
each $k$ and that for all $r \in g(k)$ and all $l > k$ there exists
$s \in g(l)$ with $D(r,s)< 2^{-k} + 2^{-l}$, since the balls of radius
$2^{-l}$ about the elements of $f(l)$ within this distance of $r$ cover
$B_{2^{-k}}(r)$. We can now define a Cauchy $\iota$-sequence $(s_k)$ by
letting $s_{k+1}$ be the $\preceq_\Uc$-minimal element of $g(k+1)$ such
that $D(s_k, s_{k+1}) < 2^{-k} + 2^{-k-1}$. Finally, define an
$\iota$-subsequence of $(r_n)$ by letting $n_{k+1}$ be the smallest index
larger than $n_k$ such that $D(s_{k+1}, r_{n_{k+1}}) < 2^{-k-1}$. This
will be a Cauchy $\iota$-sequence and hence it will converge.

(iii) $\Rightarrow$ (i): Let $z$ be an $\iota$-set of proxies of $\iota$-open
$\iota$-classes which cover $X$. Let $x$ be an $\iota$-dense
$\iota$-subset of $X$ and let $y \subseteq x \times \Qb^+$ be the
$\iota$-set of pairs $\langle p,q\rangle$ which satisfy $D(p,p') + q \leq q'$
for some $\langle p',q'\rangle$ in one of the proxies in $z$. Then the
balls $B_q(p)$ with $\langle p,q\rangle \in y$ cover $X$ and it will suffice
to find a finite subset of $y$ with the same property. Assuming $y$ is
infinite, by Theorem \ref{theorem2.10} we may write
$y = \{\langle p_n,q_n\rangle: n \in \Nb\}$ for some $\iota$-sequences
$(p_n)$ and $(q_n)$. If there exists $n$ such that every $r \in x$
satisfies $D(r, p_j) < q_j - 1/n$ for some $j \leq n$ then the balls
$\{B_{q_j}(p_j): j \leq n\}$ cover $X$ and we are done. Otherwise,
for each $n$ let $s_n$ be the $\preceq_\Uc$-least element of $x$
such that $D(s_n,p_j) \geq q_j - 1/n$ for all $j \leq n$. If (iii)
holds then we may extract a convergent $\iota$-subsequence of $(s_n)$, and
its limit $s$ must be contained in some ball $B_{q_j}(p_j)$, which yields
a contradiction. This completes the proof.
\end{proof}

\begin{defi}\label{def4.10}
An $\iota$-separable $\iota$-metric space $X$ is {\it boundedly
$\iota$-compact} if every $\iota$-closed ball $\overline{B}_a(r) =
\{s \in X: D(r,s) \leq a\}$ ($r \in X$, $a \geq 0$) is $\iota$-compact.
\end{defi}

\begin{prop}\label{prop4.11}
In any $\iota$-separable $\iota$-metric space the $\iota$-closure of any
$\iota$-set is an $\iota$-closed $\iota$-class. In any $\iota$-separable
boundedly $\iota$-compact $\iota$-metric space every $\iota$-closed
$\iota$-class is $\iota$-separable.
\end{prop}

The proof of the second part of Proposition \ref{prop4.11} is similar to
the second part of the proof of Theorem \ref{theorem3.13}, augmented by the
K\"onig's lemma technique used in the proof of (ii) $\Rightarrow$ (iii) in
Theorem \ref{theorem4.9}. (The observation that K\"onig's lemma plays an
important role in arguments of this type must be credited to the reverse
mathematics school.) The problem is to find an $\iota$-dense
$\iota$-subset $c$ of the complement of an $\iota$-open $\iota$-class $U$,
and this is done in terms of an $\iota$-dense subset $x$ of the ambient
$\iota$-metric space $X$. We construct $c$ so as to contain a point in
$\overline{B}_{2p}(r) \cap (X - U)$ for each $r \in x$ and $p \in \Qb^+$ such
that $\overline{B}_p(r) \not\subseteq U$. This is possible because bounded
$\iota$-compactness allows us to diagnose whether $\overline{B}_p(r)$ is
contained in $U$ by checking whether it is contained in finitely many
balls $B_{q'}(p')$ with $\langle p',q'\rangle$ in a proxy for $U$. Doing
this for $\overline{B}_{2^{-n}p}(s)$ for all $n \in \Nb$ and $s \in x \cap
B_p(r)$, we can then use a K\"onig's lemma argument together with the universal
well-ordering to simultaneously extract a Cauchy sequence of centers of such
balls for each $r$ and $p$ such that $\overline{B}_p(r) \not\subseteq U$. We
finally pass to the limits of the Cauchy sequences using the fact that the
map that takes convergent $\iota$-sequences to their limits is an
$\iota$-function.

\begin{defi}\label{def4.12}
Let $X$ and $Y$ be $\iota$-metric spaces. An $\iota$-function $F: X \to Y$
is {\it $\iota$-continuous} if the inverse image of every $\iota$-open
$\iota$-class in $Y$ is an $\iota$-open $\iota$-class in $X$. $F$ is
an {\it $\iota$-homeomorphism} if it is an $\iota$-bijection and both
$F$ and $F^{-1}$ are $\iota$-continuous.
\end{defi}

\begin{theo}\label{theorem4.13}
Let $X$ be an $\iota$-separable $\iota$-metric space, let $Y$ be an
$\iota$-metric space, and let $F: X \to Y$ be an $\iota$-function. Then
the following are equivalent:

(a) $F$ is $\iota$-continuous.

(b) The inverse image of every $\iota$-closed $\iota$-class in $Y$ is an
$\iota$-closed $\iota$-class in $X$.

(c) For any $\iota$-set $c \subseteq X$ with $\iota$-closure $C$, the
$\iota$-closure of $F[c]$ contains $F[C]$.

(d) $F$ preserves convergence of $\iota$-sequences.

(e) For every $r \in X$ and every $\epsilon > 0$ there exists $\delta > 0$
such that $|r - s| < \delta$ implies $|F(r) - F(s)| < \epsilon$.

If any of these conditions holds then there is an $\iota$-function which takes
any proxy $u$ of an $\iota$-open $\iota$-class $U \subseteq Y$ to a
proxy $v$ of an $\iota$-open $\iota$-class $V \subseteq X$ such that
$F^{-1}[U] = V$.
\end{theo}

\begin{prop}\label{prop4.14}
Let $X$ be an $\iota$-compact $\iota$-metric space and let $Y$ be an
$\iota$-metric space.

(a) Every $\iota$-closed $\iota$-subclass of $X$ is $\iota$-compact.

(b) If $X$ is $\iota$-separable and $F: X \to Y$ is an $\iota$-continuous
$\iota$-function then $F[X]$ is an $\iota$-separable
$\iota$-compact $\iota$-subclass of $Y$.

(c) If $X$ is $\iota$-separable and $F: X \to Y$ is an $\iota$-continuous
$\iota$-bijection then it is an $\iota$-homeomorphism.
\end{prop}

(In the proof of part (c) use parts (a) and (b) of this proposition
together with the second part of Proposition \ref{prop4.11} to show
that the image of any $\iota$-closed $\iota$-subclass of $X$ is an
$\iota$-separable $\iota$-compact $\iota$-subclass of $Y$, then show that
$Y$ must be $\iota$-separable using Theorem \ref{theorem4.13} (c), and
then use the first part of Proposition \ref{prop4.11} to verify that any
$\iota$-separable $\iota$-compact $\iota$-subclass of $Y$ must be
$\iota$-closed.)

\begin{theo}\label{theorem4.15}
(Baire category theorem) The intersection of any $\iota$-set of $\iota$-open
$\iota$-dense $\iota$-classes in an $\iota$-separable $\iota$-complete
$\iota$-metric space is $\iota$-dense.
\end{theo}

\begin{proof}
Recall that any $\iota$-set of $\iota$-open $\iota$-dense $\iota$-classes
is countable by Theorem \ref{theorem2.10}. So let $X$ be an $\iota$-separable
$\iota$-complete $\iota$-metric space and let $(u_n)$ be an $\iota$-sequence
of proxies for $\iota$-open $\iota$-dense $\iota$-classes $U_n$ in $X$. Let
$x$ be an $\iota$-dense $\iota$-subset of $X$. For each $r \in x$ and each
$q \in \Qb^+$ we can find a point $s_{r,q} \in B_q(r) \cap \bigcap U_n$ as
follows. As in the classical proof, recursively define an $\iota$-sequence
$(r_n)$ in $x$ together with an $\iota$-sequence of radii $(q_n)$ in $\Qb^+$
such that $q_n \leq 2^{-n}$ and $B_{2q_{n+1}}(r_{n+1}) \subseteq
B_{q_n}(r_n) \cap U_n$. Then let $s_{r,q}$ be the limit of this
$\iota$-sequence. Using $\preceq_\Uc$ and the fact that the map taking
convergent $\iota$-sequences to their limits is an $\iota$-function,
this construction can be carried
out simultaneously for all $r$ and $q$, so the $s_{r,q}$ constitue an
$\iota$-set which is contained in the intersection of the $U_n$ and whose
$\iota$-closure is an $\iota$-closed $\iota$-class (by Proposition
\ref{prop4.11}) which contains $x$. Since $x$ is $\iota$-dense in $X$,
we are done.
\end{proof}

\subsection{Topological spaces}\label{sect4.2}

In order to define a topology on an $\iota$-class $X$ we must specify which
of its $\iota$-subclasses are open. This is naturally done by specifying a
single $\iota$-subclass $\Tc$ of $T \times X$ where the elements of $T$
serve as proxies for the open $\iota$-classes. Thus each $a \in T$ is a
proxy for the $\iota$-class $U_a = \{r \in X: \langle a,r\rangle \in \Tc\}$.
It is convenient to simply take $T = J_2$; there is no loss in generality
since we always have $T \times X \subseteq J_2 \times X$.

Another feature of the following definition might require explanation.
Classically the family of open sets is
closed under arbitrary unions and finite intersections. However, verifying
these closure properties in any given example would generally be done by in
effect determining proxies for the union and intersection. Thus it is natural
in our setting, and it bars no important examples, to require the existence
of $\iota$-functions which evaluate unions and intersections of proxies.

\begin{defi}\label{def4.16}
An {\it $\iota$-topological space} (or just {\it $\iota$-space}) is an
$\iota$-class $X$ together with an $\iota$-subclass $\Tc$ of
$J_2 \times X$ with the following properties:

(i) $\emptyset = U_a$ and $X = U_b$ for some $a,b \in J_2$;

(ii) there is an $\iota$-function $\kappa: \Pc_\iota(J_2) \to J_2$ such
that for any $\iota$-set $x$ we have
$$\bigcup_{a \in x} U_a = U_{\kappa(x)};$$

(iii) there is an $\iota$-function $\lambda: J_2 \times J_2 \to J_2$ such
that for any $\iota$-sets $a$ and $b$ we have
$$U_a \cap U_b = U_{\lambda(a,b)}$$

\noindent where $U_a = \{r \in X: \langle a,r\rangle \in \Tc\}$ for any
$a \in J_2$. The $U_a$ are the {\it $\iota$-open $\iota$-classes} in
$X$ and $\Tc$ is an {\it $\iota$-topology} on $X$. The {\it $\iota$-closed
$\iota$-classes} are the complements $X - U_a$.
\end{defi}

\begin{exam}\label{example4.17}
(a) Let $x$ be an $\iota$-set and let $\Tc$ be the $\iota$-class of all
ordered pairs $\langle y,r\rangle \in J_2 \times x$ such that
$r \in x \cap y$. This is the {\rm discrete $\iota$-topology}
on $x$.

(b) Let $X$ be an $\iota$-class and let $\Tc$ be the $\iota$-class of all
ordered pairs $\langle x,r\rangle \in J_2 \times X$ such that $x =
\langle 0,y\rangle$ for some $y \in J_2$ and $r \in X - y$. This is the
{\rm co-countable $\iota$-topology} on $X$.
\end{exam}

It is straightforward to verify that the discrete and co-countable
$\iota$-topologies satisfy Definition \ref{def4.16}.

\begin{defi}\label{def4.18}
Let $X$ be an $\iota$-space with $\iota$-topology
$\Tc \subseteq J_2 \times X$ and let $Y$ be an $\iota$-subclass
of $X$. Then $\Tc' = \Tc \cap (J_2 \times Y)$ is the {\it relative
$\iota$-topology} on $Y$ and $Y$ equipped with this $\iota$-topology is
an {\it $\iota$-subspace} of $X$.
\end{defi}

Next we introduce a basic tool for constructing $\iota$-topologies.
If $\Bc \subseteq B \times X$ and $a \in B$ then we
write $B_a = \{r \in X: \langle a,r\rangle \in \Bc\}$.

\begin{prop}\label{prop4.19}
Let $X$ and $B$ be $\iota$-classes and let $\Bc$ be an $\iota$-subclass of
$B \times X$ such that $B_a = X$ for some $a \in B$. Suppose also that
there exists an $\iota$-function $\lambda_0: B \times B \to \Pc_\iota(B)$
such that for any $a,b \in B$ we have $B_a \cap B_b
= \bigcup_{c \in \lambda_0(a,b)} B_c$. Then the $\iota$-class $\Tc
\subseteq J_2 \times X$ consisting of the ordered pairs $\langle x,r\rangle$
such that $x \in \Pc_\iota(B)$ and $r \in \bigcup_{a \in x} B_a$ is an
$\iota$-topology on $X$.
\end{prop}

\begin{defi}\label{def4.20}
The $\iota$-topology defined in Proposition \ref{prop4.19} is the
$\iota$-topology {\it generated by $\Bc$}.
\end{defi}

\begin{exam}\label{example4.21}
Let $X$ be an $\iota$-separable $\iota$-metric space and let $x \subseteq X$
be an $\iota$-dense $\iota$-subset. Let $B = x \times \Qb^+$ and define
$\Bc \subseteq B \times X$ to be $\{\langle r,q, s\rangle: D(r,s) < q\}$.
This satisfies the conditions of Proposition \ref{prop4.19} with
$$\lambda_0(\{\langle r_1, q_1\rangle, \langle r_2, q_2\rangle\})
= \{\langle r,q\rangle \in x \times \Qb^+:
D(r, r_i) + q \leq q_i\hbox{ for }i = 1,2\}.$$
The $\iota$-open classes for the $\iota$-topology generated by $\Bc$ are
precisely the $\iota$-open classes identified in Definition \ref{def4.5}.
This is the {\rm metric $\iota$-topology} on $X$.
\end{exam}

\begin{defi}\label{def4.22}
(a) An {\it $\iota$-family of $\iota$-topological spaces} consists of
an $\iota$-set $x$, an $\iota$-class $\Xc \subseteq x \times J_2$, and
an $\iota$-class $\Tc \subseteq x \times J_2 \times J_2$ with the
property that for each $a \in x$ the $\iota$-class
$$\Tc_a = \{\langle b, r\rangle: \langle a,b,r\rangle \in \Tc\}$$
is an $\iota$-topology on $X_a = \{r \in J_2: \langle a,r\rangle \in \Xc\}$.
We will write $\Xc = \{X_a: a \in x\}$.

(b) Given an $\iota$-family of $\iota$-topological spaces $\{X_a: a \in x\}$,
let $B$ be the $\iota$-class of functions from finite subsets of $x$ into
$J_2$ and let $\Bc \subseteq B \times {\prod}_\iota X_a$ be the
$\iota$-class of pairs $\langle h,f\rangle$ such that for each
$a \in {\rm dom}(h)$ we have $\langle a, h(a), f(a)\rangle \in \Tc$.
The {\it product $\iota$-topology} on ${\prod}_\iota X_a$ is the
$\iota$-topology generated by $\Bc$.
\end{defi}

\subsection{Continuity}\label{sect4.3}

\begin{defi}\label{def4.23}
Let $X$ and $Y$ be $\iota$-topological spaces with $\iota$-topologies $\Tc_X$
and $\Tc_Y$. An $\iota$-function $F: X \to Y$ is {\it $\iota$-continuous} if
there is an $\iota$-function $\widetilde{F}: J_2 \to J_2$ such that for
each $b \in J_2$ we have $F^{-1}[V_b] = U_{\widetilde{F}(b)}$, where
$U_a = \{r \in X: \langle a,r\rangle \in \Tc_X\}$ and
$V_b = \{s \in Y: \langle b,s\rangle \in \Tc_Y\}$.
An {\it $\iota$-homeomorphism} is an $\iota$-continuous $\iota$-bijection
whose inverse is also $\iota$-continuous.
\end{defi}

\begin{prop}\label{prop4.24}
Compositions of $\iota$-continuous $\iota$-functions are $\iota$-continuous.
\end{prop}

\begin{defi}\label{def4.25}
Let $X$ be an $\iota$-space with $\iota$-topology $\Tc \subseteq J_2 \times X$
and let $B$ be an $\iota$-class. Then $\Bc = \Tc \cap (B \times X)$ is an
{\it $\iota$-base} for $\Tc$ if it
satisfies the conditions of Proposition \ref{prop4.19} and the
identity map on $X$ is an $\iota$-homeomorphism between $\Tc$ and the
$\iota$-topology generated by $\Bc$. $X$ is {\it $\iota$-second countable}
if it has an $\iota$-base for which $B$ is an $\iota$-set.
\end{defi}

\begin{prop}\label{prop4.26}
Let $X$ and $Y$ be $\iota$-topological spaces and let $F: X \to Y$ be
an $\iota$-function. Suppose $\Bc \subseteq B \times Y$ is an $\iota$-base
for $Y$. Then $F$ is $\iota$-continuous if and only if there is an
$\iota$-function $\widetilde{F}_0: B \to J_2$ such that for each $b \in B$
we have $F^{-1}[V_b] = U_{\widetilde{F}_0(b)}$, where
$U_a = \{r \in X: \langle a,r\rangle \in \Tc_X\}$ and
$V_b = \{s \in Y: \langle b,s\rangle \in \Tc_Y\}$.
\end{prop}

The key observation in the proof of Proposition \ref{prop4.26} is that
if $\Bc$ generates $\Tc_Y$ (it is sufficient to consider this case by
Proposition \ref{prop4.24}) then the function from $\Pc_\iota(B)$ to $J_2$
which takes an $\iota$-subset $x \subseteq B$ to $\{\widetilde{F}(b):
b \in x\}$ is an $\iota$-function. We can then define an $\iota$-function
$\widetilde{F}$ that verifies $\iota$-continuity by composing this function
with the function $\kappa$ of Definition \ref{def4.16} (relative to $\Tc_X$).

\begin{prop}\label{prop4.27}
Let $\{X_a: a \in x\}$ be an $\iota$-family of $\iota$-topological
spaces and let $Y$ be an $\iota$-topological space. For each $b \in x$ let
$\pi_b: {\prod}_\iota X_a \to X_b$ be the natural projection.

(a) Each $\pi_b$ is $\iota$-continuous.

(b) Let $F: Y \to {\prod}_\iota X_a$ be an $\iota$-function. Then $F$ is
$\iota$-continuous if and only if the $\iota$-function
$G: Y \times x \to \coprod X_a$ (= the disjoint union of the $X_a$)
defined by $G(r,a) = F(r)(a) \in X_a$ is $\iota$-continuous. Here we give $x$
the discrete $\iota$-topology, $Y \times x$ the product $\iota$-topology, and
$\coprod X_a$ the $\iota$-topology generated by the $\iota$-topologies on the
individual $\iota$-spaces.
\end{prop}

At this point we could go on to develop general topology in the $J_2$
setting. Most standard results go through, although usually additional
hypotheses such as separability or second countability are required.
However, this seems somewhat extraneous to the development of core
mathematics since most topological spaces that appear in ordinary
settings are separable and metrizable and for these spaces most of the
basic results in general topology are easy. Those spaces not of this type
which do appear in mainstream settings (e.g., the Zariski topology)
do not seem to require deep results from general topology. The
weak* topology on the dual of a separable Banach space is typically
not metrizable but its restriction to the unit ball is, and this coupled
with the Krein-Smullian theorem appears to render metric space theory
sufficient for most applications.

\section{Other topics}\label{sect5}

\begin{quote}
{\it We sketch ways of developing various other topics from abstract
analysis within $J_2$.}
\end{quote}

\subsection{Measure and integration}\label{sect5.1}

We can define {\it $\iota$-$\sigma$-algebras} in a manner analogous to the
definition of $\iota$-topologies (Definition \ref{def4.16}). Thus, an
$\iota$-$\sigma$-algebra on an $\iota$-class $X$ is an $\iota$-subclass $\MC$
of $J_2 \times X$ for which there exist $\iota$-functions $\kappa: J_2 \to J_2$
and $\lambda: \Pc_\iota(J_2) \to J_2$ such that for any $a$ and $x$ we have
$X - M_a = M_{\kappa(a)}$ and $\bigcup_{b \in x} M_b = M_{\lambda(x)}$, where
$M_a = \{r \in X: \langle a,r\rangle \in \MC\}$ for any $a \in J_2$.
We call the sets $M_a$ the {\it measurable $\iota$-subclasses of $X$}
and we call $a$ a {\it proxy} for $M_a$.

Most $\sigma$-algebras of interest are defined in terms of a generating
set, and it may not be obvious how to generate an $\iota$-$\sigma$-algebra
from a given family of $\iota$-classes because classically this involves a
recursive construction along the set of all countable ordinals, which are
not available in $J_2$. However, it can be done in the following way.

Let $X$ and $B$ be $\iota$-classes and let $\Bc \subseteq B \times X$ be an
$\iota$-class. We describe a way to construct an $\iota$-$\sigma$-algebra
$\MC$ generated by the sets $B_a = \{r \in X: \langle a,r\rangle \in \Bc\}$.
We may assume that the set $\{B_a: a \in B\}$ is an algebra, i.e., it is
closed under complements and finite unions. It also simplifies matters
slightly to observe that we only need closure under the single operation
which, for any $\iota$-set of $\iota$-measurable $\iota$-subclasses of $X$,
forms the complement of their union. Applying this operation to an
$\iota$-set containing only one $\iota$-class yields the complement of that
$\iota$-class, and composing the operation with complementation then produces
unions. Thus, proxies for $\iota$-measurable $\iota$-classes will be given by
trees whose terminal nodes are elements of $B$ (i.e., proxies for generating
$\iota$-classes) and each of whose non-terminal nodes will represent the
complement of the union of the $\iota$-classes corresponding to its
immediate successors. Specifically, define a {\it $B$-tree}
to be a small $\iota$-function $f$ whose range is contained in $\Nb$ and
which satisfies the following conditions. For all $a \in {\rm dom}(f)$, if
$f(a) = 0$ then $a \in B$; if $f(a) = 1$ then $a$ is infinite; and if
$f(a) \geq 1$ then $a \subseteq \Gamma(f)$ and
we have $f(a) = 1 + \sup_{b \in {\rm dom}(a)} f(b)$. We also require that
there exist exactly one element on which $f$ attains a maximal value.

We now define $\MC \subseteq J_2 \times X$. If $a \in J_2$ is not a
$B$-tree then let $M_a = \emptyset$. For any $B$-tree $f$, we define
$M_f$ to be the set of $r \in X$ for which there exists an $\iota$-function
$g: {\rm dom}(f) \to \{0,1\}$ such that (1) if $f(a) = 0$ then
$g(a) = 1$ if and only if $r \in B_a$, (2) if $f(a) > 1$ then $g(a) = 1$ if
and only if $g(b) = 0$ for all $b \in {\rm dom}(a)$, and (3) $g(a) = 1$ for
the element $a$ on which $f$ attains its maximum value.

What makes this construction work is the condition that $a$ be infinite
if $f(a) = 1$. This means that either the domain of $f$ consists of a single
element in $B$ or else, in set-theoretic terms, $f$ has infinite rank. An
easy induction shows that every set in $J_2$ has rank less than $2\omega$,
so it follows that any $\iota$-set of $B$-trees contains only
$B$-trees of at most some maximal finite height. Thus, any $\iota$-set of
$B$-trees can be amalgamated into a single $B$-tree and this can
be used to show that $\MC$ is closed under complements of unions.

We define an {\it $\iota$-measure} on an $\iota$-class $X$ equipped with an
$\iota$-$\sigma$-algebra $\MC$ to be an $\iota$-function $\mu$ from $J_2$ into
$[0,\infty]$ which satisfies $\mu(a) = \mu(b)$ if $M_a = M_b$, $\mu(a) = 0$
if $M_a = \emptyset$, and $\mu(\kappa(x)) = \sum_{a \in x} \mu(a)$ if
the $\iota$-classes $M_a$ ($a \in x$) are disjoint. If $\MC$ is generated
by an algebra $\Bc$ and $\mu_0$ is a premeasure on $\Bc$ then we can extend
$\mu_0$ to an $\iota$-measure $\mu$ by a standard inner/outer measure
construction. Specifically, if $\mu_0(X) < \infty$ then we can show by
induction on $B$-trees that for every $M_a$ there exists a pair of
$\iota$-sequences $(a_n)$ and $(b_n)$ such that each $a_n$ is an
$\iota$-subset of
$B$ with $M_a \subseteq \bigcup_{b \in a_n} B_b$; each $b_n$ is an
$\iota$-subset of $B$ with $X - M_a \subseteq \bigcup_{b \in b_n} B_b$;
and $\sum_{b \in a_n} \mu_0(b) + \sum_{b \in b_n} \mu_0(b)
\leq \mu_0(X) + 2^{-n}$. This can
then be used to define $\mu$. The case $\mu_0(X) = \infty$ introduces no
fundamental difficulties.

The theory of integration seems to carry over with no major complications.
We define an {\it $\iota$-measurable $\iota$-function} from $X$ to $Y$ to be
an $\iota$-function $F: X \to Y$ such that there exists an $\iota$-function
$\widetilde{F}$ from proxies of $\iota$-measurable $\iota$-classes $N_b$ in $Y$
to proxies of $\iota$-measurable $\iota$-classes $M_a$ in $X$ which satisfies
$F^{-1}[N_b] = M_{\widetilde{F}(b)}$ for all $b \in J_2$. This makes
approximation by simple functions possible since, e.g., $\{(-\infty, k/n]
\cap \Rb_\iota: k \in \Zb\}$ is an $\iota$-set of $\iota$-measurable
$\iota$-classes in
$\Rb_\iota$ for each $n$, and the family of all such $\iota$-sets is an
$\iota$-set, so the fact that $\widetilde{F}$ is an $\iota$-function would
permit the construction of a corresponding $\iota$-set of simple
$\iota$-functions which approximate $F$ and have range contained in
$\{k/n: k \in \Zb\}$.

\subsection{Banach spaces}\label{sect5.2}
We define an {\it $\iota$-real $\iota$-vector space} to be
an $\iota$-class $V$ equipped with $\iota$-functions
$+: V \times V \to V$ and $\cdot: \Rb_\iota \times V \to V$ which satisfy
the usual vector space axioms, such that for any $\iota$-set $x \subseteq V$
and any small $\iota$-subfield $y \subseteq \Rb_\iota$, the $y$-linear span
of $x$ is an $\iota$-set. An {\it $\iota$-real $\iota$-Banach space} is an
$\iota$-real $\iota$-vector space equipped with
an $\iota$-function $\|\cdot\|: V \to
[0,\infty)$ which satisfies the norm axioms and such that $D(r,s) = \|r - s\|$
is an $\iota$-complete $\iota$-metric on $V$. We can define $\iota$-complex
$\iota$-Banach spaces analogously.

Much of the general theory of Banach spaces seems to carry over fairly
easily. There is no problem proving versions of the open mapping and
closed graph theorems and the principle of uniform boundedness.
The Hahn-Banach theorem presents a difficulty, however. Suppose
we want to extend a bounded linear functional on an $\iota$-separable
$\iota$-subspace of an $\iota$-separable $\iota$-real $\iota$-Banach space $V$ to the
entire space without increasing its norm. Extending by a single dimension
can be done in the usual way, but it is not obvious that a sequence of
one-dimensional extensions will give rise to a (small) $\iota$-function on
an $\iota$-dense $\iota$-subset of $V$. The problem is that at each step
we define a new $\iota$-real number and it is not clear that the sequence
of $\iota$-reals so obtained must be an $\iota$-set. One way to get around
this problem is by using only rational numbers to define the extensions and
allowing the norm to increase by a small amount. In fact, if we use
rationals for the extensions we can simultaneously, for all $n$, define
extensions which increase the norm by at most $2^{-n}$ and which converge
pointwise to a single extension that does not increase the norm. In this
way we can ensure that the extension to $V$ is an $\iota$-function.

If $V$ is an $\iota$-separable $\iota$-Banach space then we define its
{\it $\iota$-dual} to be the $\iota$-class of bounded linear $\iota$-functions
from
an $\iota$-dense $\iota$-subset of $V$ (which is without loss of generality a
vector space over $\Qb$, say) into the scalar field. Each such $\iota$-function
$f$ can be regarded as a proxy for a bounded linear $\iota$-function $F$
from $V$ into the scalars.

\subsection{Function spaces}\label{sect5.3}

Let $X$ be an $\iota$-separable
$\iota$-compact $\iota$-metric space. We can define $C(X)$
to be the $\iota$-class of uniformly $\iota$-continuous $\iota$-functions
from an $\iota$-dense $\iota$-subset of $X$ into the scalars. We
want to think of each such function as a proxy for a continuous function
on $X$ and it is clear how to do this. The following seems like a good
general definition. An {\it $\iota$-function space} over $X$ is an
$\iota$-subclass $\Fc$ of $V \times X \times \Fb_\iota$, where $V$ is
an $\iota$-Banach space, $\Fb_\iota$ is the scalar field, for each
$v \in V$ the $\iota$-class $\{\langle r, s\rangle \in X \times \Fb_\iota:
\langle v,r,s\rangle \in \Fc\}$ is the graph of an $\iota$-function
$F_v: X \to \Fb_\iota$, and we have $F_{v+w} = F_v + F_w$ and
$F_{rv} = rF_v$ for all $v,w \in V$ and $r \in \Fb_\iota$. That is,
the vector space operations in $V$ correspond to pointwise operations on
the $\iota$-functions $F_v$.

This definition encompasses not only standard function spaces like
$C(X)$, but also the $\iota$-dual of an $\iota$-Banach space as defined
in the preceding section. However, it does not work with ``function''
spaces like $L^p(X)$ whose elements are actually equivalence classes of
functions. One way to handle these would be to weaken the $\iota$-function
space definition so that $F_{v+w} = F_v + F_w$ and $F_{rv} = rF_v$ hold
off of an $\iota$-class with $\iota$-measure zero, for each $v,w,r$.

Versions of the Stone-Weierstrass and Tietze extension theorems hold for
$C(X)$ with $X$ an $\iota$-separable
$\iota$-compact $\iota$-metric space. In the
former case, the hypothesis must include that the $\iota$-dense subalgebra
of $C(X)$ be $\iota$-separable. We can then work exclusively with an
$\iota$-dense $\iota$-set of functions in the subalgebra, which permits
execution of the compactness arguments needed in the proof.

%##

\end{document}